\documentclass[11pt,twoside,A4paper]{article} 
\usepackage{amssymb,amsthm,amsmath,amscd,amsfonts,txfonts}
\usepackage{xcolor}

\usepackage[tmargin=1.5cm,bmargin=1.5cm,lmargin=1.7cm,rmargin=1.7cm]{geometry}

\theoremstyle{plain}
\newtheorem{lemma}{Lemma}[section]
\newtheorem{proposition}{Proposition}[section]
\newtheorem{theorem}{Theorem}[section]
\newtheorem{corollary}{Corollary}[section]

\theoremstyle{definition}
\newtheorem{definition}{Definition}[section]

\newtheorem{example}{Example}[section]

\theoremstyle{remark}
\newtheorem{remark}{Remark}[section]

\title{\bfseries\scshape{Some results on compatible ternary Leibniz algebras}}
\author{\bfseries\scshape Kol B\'eatrice GAMOU\thanks{E-mail address: \tt{kolbeatrice18@gmail.com}} \\
D\'epartement de Math\'ematiques,\\
Universit\'e Gamal Abdel Nasser, Conakry, Guinea. \\
\bfseries\scshape  Ibrahima BAKAYOKO \thanks{E-mail address: \tt{ibrahimabakayoko27@gmail.com}}\\
D\'epartement de Math\'ematiques,\\ 
Universit\'e de N'Z\'er\'ekor\'e, Guinea.}

\date{} 

\begin{document} 
\maketitle

\tableofcontents

\begin{abstract} 
In this paper, we introduce compatible ternary Leibniz algebras, (dual)Nijenhuis pairs from the second-order deformation of ternary Leibniz
 algebras with a representarion and study the  invariance of certains operators (generalized derivations, Rota-Baxter operators, 
Reynolds operators, element of centroid, averaging operators, Nijenhuis operators) whenever we go from 
Leibniz algebras to ternary Leibniz algebras. We also give contructions of (compatible) ternary Leibniz algebras either from 
averaging operators, modules or $O$-operator.

\end{abstract} 


\section{Introduction}
A ternary algebra consists of a linear vector space $V$ together with a trilinear map $\mu : V\times V\times V\rightarrow V$.
 In other words, ternary algebras are vector spaces equipped with a multiplication with three  items instead of two, as in classical algebraic 
structures. They originate from the work of Jacobson in 1949 in the study of associative algebra $(A, .)$ that
 are closed relative to the ternary operation $[[a, b], c]$, where $[a, b]=ab-ba$. According to some conditions satisfyied by the multiplication 
$\mu$, we dispose of  (partial or total) ternary associative algebras, ternary Leibniz algebras, $3$-Lie algebras, ternary Leibniz-Poisson 
algebras, ternary Hopf algebras, ternary Heap algebras,  Comstrans algebras, Akivis algebras, 
Lie-Yamaguti algebras, Lie triple systems \cite{HB, HT1, TH1}, ternary Jordan algebras \cite{IAP}, Jordan-Lie triple systems \cite{SN},
 Jordan triple and so on.

It is known that certains algebraic structures admit left and right version such as right symmetric algebras (left symmetric algebras), right
Leibniz algebras (left Leibniz algebras), right BiHom-Lie algebras (right BiHom-Lie algebras) and so on. 
In this paper, we deal with right $3$-Leibniz algebras, also called right ternary  Leibniz algebras and simply refer it ternary 
Leibniz algebra in this paper. It reads
\begin{eqnarray}
 [[x, y, z], t, u]=[x, y, [z, t, u]]+[x, [y, t, u], z]+[[x, t, u], y, z],\nonumber
\end{eqnarray}
for any $x, y, z\in V$. Whenever the bracket $[-, -, -]$ is skew-symmetric with respect to any pair of variables, $(V, [-, -, -])$ is said to be
  a ternary Lie algebra or $3$-Lie \cite{AKM, JAS, BW, BRB, SJ}. 

It is well-known that mathematical objects are often understood through studying operators defined on them. For instance, in Gallois theory a field
 is studied by its automorphisms, in analysis functions are studied through their derivations,  and in geometry manifolds are studied through 
their vector fields. Fifty years ago, several operators have been found from studies in analysis, probability and combinatorics. Among these
 operators, one can cite, element of centroid \cite{RM2}, averaging operator, Reynolds operator, Leroux's TD operator, Nijenhuis operator and
 Rota-Baxter operator \cite{LG, MD}.

The Rota-Baxter operator originated from the work of G. Baxter \cite{GB, AM} on Spitzer's identity\cite{F} in fluctuation theory.
Rota- Baxter algebras (associative algebra with Rota-Baxter operator) are used in many fiels of mathematics and mathematical Physics. 
In mathematics, they are used in algebra,
 number theory, operads and combinatorics \cite{MA}, \cite{CA1}, \cite{CA2}, \cite{PC}, \cite{GCB}. In 
mathematical physics they appear as the operator form of the classical Yang Baxter equation \cite{CA1} or as the fondamental algebraic strucutre in
 the normalisation of quantum field theory of Connes and Kreimer \cite{CK}. In non-associative algebra, the Rota-Baxter operators are used  in order 
to produce another one of the same type or not from the previous one.

The Nijenhuis operator on an associative algebra was introduced in \cite{CJ} in the study quantum bi-Hamiltonian
systems while the notion of Nijenhuis operator on a Lie algebra originated from the concept of Nijenhuis tensor that was introduced by Nijenhuis
in the study of pseudo-complex manifolds and was related to the well known concepts of Schouten-Nijenhuis bracket, the Frolicher-Nijenhuis
 bracket \cite{FN}, and the Nijenhuis-Richardson bracket. The associative analog of the Nijenhuis relation may be regaded as the homogeneous version of 
Rota-Baxter relation\cite{PL}.

Leibniz algebras, also known as Loday algebras, are introduced by Loday in \cite{JL} as a non-commutative version of Lie algebras.
Various aspects of these algebras are intensively investigated these last dacades. For instance, Almutari H. and  Ahmad A. G. 
studied some properties of centroid and quasicentroids \cite{HA2} for Leibniz algebras and determine the centroids of low-dimensional Leibniz
algebras and give their classification. While Hassan A. and  ABD G. A.\cite{}  deal with some properties and classification of complexe Leibniz
algebras of dimensions less or equal to 4.  In \cite{MS}
the authors characterize compatible Leibniz algebras in terms of Maurer-Cartan elements of a suitable differential graded
Lie algebra, define a cohomology theory of compatible Leibniz algebras and study
the abelian extension of compatible Leibniz algebras.

The aim of this paper is to study invariance of operators when we pass from Leibniz algebra to the associated ternary Leibniz algebra,
the compatibility of ternary Leibniz algebras and their modules, the infinitesimal deformation of ternary Leibniz algebras with representations.
First, we introduce compatible ternary Leibniz algebras, gives some constructions using averaging operators and establish some properties
connecting compatible Leibniz algebra to the associated ternary Leibniz algebra. 
Then, we establish that operators on a given Leibniz algebras keep his name on the associated ternary Leibniz algebras. For istance, a Reynolds 
operator on a Leibniz algebra is also a Reynolds operator on the associated ternary Leibniz algebra. Some examples of operators are given and 
the notion of $O$-operators are introduced for ternary Leibniz algebras and a new contruction of ternary Leibniz algebras from $O$-operators is done.
Next, we study modules and representations of (compatible) ternary Leibniz algebras as well the compatibility concept of modules and
 representations of compatible ternary Leibniz algebras. Finally,we introduce Nijenhuis operators and (dual)-Nijenhuis pair for ternary Leibniz 
algebras.
 
The paper is organized as follows. In section 2, we recall basics notions on Leibniz algebras and ternary Leibniz algebras i.e. definitions,
examples, some construction methods, modules and representation. 
Section 3 concerns the constructions of (compatible) ternary Leibniz algebra from averaging operators, the preservation of the operators from
Leibniz algebras to the associated ternary Leibniz algebras.
 The section 4 is devoted to the modules and the representations of compatible ternary Leibniz algebras and the concept of Nijenhuis pair on 
ternary Leibniz algebras.
\section{Preliminaries}
This section is devoted to the recall of some results on (compatible) Leibniz algebras and their modules.
\subsection{Leibniz  algebras}
This subsection is devoted to the recall of some results on (compatible) Leibniz algebras and their modules.
\begin{definition}\label{}
Let $L$ be a vector vector over a field $\mathbb{K}$  space and  $x, y, z\in  L$.
 a) A Leibniz algebra structure on $L$ is a  bilinear map 
$[-, -] : L\otimes L\rightarrow L$   satisfying
 \begin{eqnarray}
   [[x, y], z]=[x, [y, z]]+[[x, z], y] \label{cpa}
 \end{eqnarray}
b)  Two Leibniz algebras $(L, [-, -]_1)$ and $(L, [-, -]_2)$   are called compatible if for any $\lambda_1, \lambda_2\in\mathbb{K}$,  
the following bracket 
\begin{eqnarray}
[x, y]=\lambda_1[x, y]_1+\lambda_2[x, y]_2,  \: \: \forall x, y\in L,\label{cla1}
\end{eqnarray}
  defines a Leibniz algebra structure on $L$.

\end{definition}
See \cite{MS} for examples.

\begin{remark}
 The braket $(\ref{cla1})$ defines a Leibniz algebra structure on $L$ if and only if
 \begin{eqnarray}
[[x, y]_1, z]_2+ [[x, y]_2, z]_1 =[x, [y, z]_1]_2+[x, [y, z]_2]_1 +[[x, z]_1, y]_2+[[x, z]_2, y]_1. \label{cla2}
\end{eqnarray}
\end{remark}

\begin{definition}\label{}
The triple $(L, [-, -]_1, [-, -]_2)$ is said to be a compatible Leibniz algebra if $(L, [-, -]_1)$
and $(L, [-, -]_2)$ are both Leibniz algebras and (\ref{cla2}) holds.
\end{definition}

\begin{definition}
 Let $(L, [-, -])$ be a Leibniz algebra. A linear map $\beta : L\rightarrow L$ is said to be an averaging operator if 
 $$\beta([\beta(x), y])=[\beta(x), \beta(y)]=\beta([x, \beta(y)]),$$
for any $x, y\in L$.
\end{definition}

The following definition is inspired from (\cite{AD}).
\begin{definition}
Two averaging operators on a Hom-Leibniz algebras $L$ are said to be compatible if their sum is also an averaging operator on $L$.
This means that 
 \begin{eqnarray}
  [\beta_2([\beta_1(x), y]), z]+[\beta_1([\beta_2(x), y]), z]=[[\beta_1(x), \beta_2(y)], z]+[[\beta_2(x), \beta_1(y)], z] \label{cao}
 \end{eqnarray}
for all $x, y, z\in L$.
\end{definition}

We have already proved that $[-, -]_\beta$ is a Leibniz algebra (\cite{IS}, Proposition 3). Now we have the below proposition.
\begin{proposition}\label{}
Let $(L, [-, -])$ be a Leibniz algebra and $\beta_1, \beta_2 :L\rightarrow L$ two compatible commuting injective averaging 
operators on $L$. Then, $(L, [-, -]_{\beta_1}, [-, -]_{\beta_2})$ is a compatible Leibniz algebra.
\end{proposition}
\begin{proof}
The compatibility comes from
 a direct computation.
\end{proof}

\begin{definition}
 Let $(L, [-, -])$ be a Leibniz algebra. A linear map $N : L\rightarrow L$ is said to be 
a Nijenhuis operator if 
$$[N(x), N(y)]=N\Big([N(x), y]+[x, N(y)]-N([x, y])\Big),$$
for any $x, y\in L$.
\end{definition}

\begin{proposition}\label{}
Let $(L, [-, -])$ be a Leibniz algebra. Then, $(L, [-, -], [-, -]_N)$ is a compatible Leibniz algebra, where for any $x, y\in L$,
$$
[x, y]_N=[N(x), y]+[x, N(y)]+[x, y].
$$
\end{proposition}
\begin{proof}
 It is clear from (\cite{IS}, Proposition 7) that $(L, [-, -]_N)$ is a Leibniz algebra. The compatibility comes from
 a straightforward computation.
\end{proof}

Now we recall modules over (compatible) Leibniz algebras (\cite{MS}).
\begin{definition}
 Let $(L, [-, -])$ be a Leibniz algebra,  $M$ be a linear vector space, $l : L\times M\rightarrow M$ and
$r : M\times L\rightarrow M$ be two bilinear maps. The triple $(M, l, r)$ is said to be a bimodule over $L$ if 
\begin{eqnarray}
  l([x, y], m) &=& l(x, l(y, m))+r(l(x, m), y),  \cr 
 r(l(x, m), y)  &=& l(x, r(m, y))+l([x, y], m),\cr
 r(r(m, x), y)&=&r(m, [x, y])+r(r(m, y), x),
\end{eqnarray}
for all $x, y\in L, m\in M$.
\end{definition}

\begin{definition}
 Let $(L, [-, -])$ be a Leibniz algebra, $R : L\rightarrow L$ a Rota-Baxter operator on $L$, $(M, l, r)$ a bimodule on $M$ 
 and $R_M : M\rightarrow M$ a linear map on $M$. We say that $(M, l, r, R_M)$ is a Rota-Baxter bimodule on $L$ if, for any $x, y\in L$,
$m\in M$, we have 
\begin{eqnarray}
 l(R(x), R_M(m))&=&R_M\Big(l(R(x), m)+l(x, R_M(m))\Big),\nonumber\\
r(R_M(m), R(x))&=&R_M\Big(r(R_M(m), x)+r(m, R(x))\Big).\nonumber
\end{eqnarray}
\end{definition}

The proof of the following proposition is straightforward.
\begin{proposition}
 Let $(M, l, r, R_M)$ be a Rota-Baxter bimodule over the Rota-Baxter Leibniz algebra $(L, [-, -], R)$.
Let us define the two bilinear maps $l' : L\times M\rightarrow, r' : M\times L\rightarrow M$ as follows
 \begin{eqnarray}
l'(x, m)&=&l(R(x), m)+l(x, R_M(m)),\nonumber\\
 r'(m, x)&=&r(R_M(m), x)+r(m, R(x)),\nonumber
\end{eqnarray}
for any $x\in L, m\in M$. Then $(M, l', r')$ is a bimodule over the Leibniz algebra $L_R=(L, [-, -]_R)$.
\end{proposition}

\begin{definition}
 Let $(L, [-, -]_1, [-, -]_2)$ be a compatible Leibniz algebra. A compatible bimodule is a a linear vector space $M$
 endowed with four even linear maps 
\begin{eqnarray}
\begin{array}{lll}
l_1  : L\times M\rightarrow M, \: (x, m)\mapsto l_1(x, m) & ; &   r_1: M\times L\rightarrow M, \: (m, x)\mapsto r_1(m, x)\\
 l_2 : L\times M\rightarrow M, \: (x, m)\mapsto l_2(x, m) & ; &  r_2 : M\times L\rightarrow M, \: (m, x)\mapsto r_2(m, x),
\end{array}
\end{eqnarray}
such that $(M, l_1, r_1)$ is a bimodule over $(L, [-, -]_1)$, $(M, l_2, r_2)$ is a bimodule over $(L, [-, -]_2)$ and
  for any $x, y\in L$, $m\in M$,
\begin{eqnarray}
 l_2([x, y]_1, m)+ l_1([x, y]_2, m) &=& l_2(x, l_1(y, m))+l_1(x, l_2(y, m))+r_2(l_1(x, m), y)+r_1(r_2(x, m), y),  \cr 
 r_2(l_1(x, m), y)+ r_1(l_2(x, m), y)  &=& l_2(x, r_1(m, y))+l_2(x, r_2(m, y)) +l_2([x, y]_1, m)+l_1([x, y]_2, m),\cr
r_2(r_1(m, x), y)+ r_1(r_2(m, x), y)&=&r_2(m, [x, y]_1)+r_1(m, [x, y]_2)
 +r_2(r_1(m, y), x)+r_1(r_2(m, y), x).\nonumber
\end{eqnarray}
\end{definition}

\begin{example}
 Any compatible Leibniz algebra is a compatible bimodule over itself.
\end{example}
The following theorem comes  from a straightforward computation.
\begin{theorem}
 Let $(M, l_1, r_1, l_2, r_2)$ be a compatible bimodule over the compatible Leibniz algebra
$(L, [-, -]_1, [-, -]_2)$. Then, $L\oplus M$ has a compatible Leibniz structure with the following operations
\begin{eqnarray}
 [x+m_1, y+m_2]_i:=[x, y]_i+l_i(x, m_2)+r_i(m_1, y),
\end{eqnarray}
for all $x, y\in L, m_i\in M, i=1, 2$.
\end{theorem}

The below definition will be useful for Theorem \ref{ops}.
\begin{definition}
 Let $(L, [-, -])$ be a Leibniz algebra. A linear map $\beta : L\rightarrow L$ is said to be 
\begin{enumerate}
 \item [1)] a {\bf derivation} of weight $\lambda\in\mathbb{K}$ if \hspace{0,25cm} 
$\beta([x, y])=[\beta(x), y]+[x, \beta(y)]+\lambda[x, y],$
\item [2)] a {\bf Rota-Baxter operator of weight $\lambda\in\mathbb{K}$} if \hspace{0.25cm}
 $[\beta(x), \beta(y)]=\beta\Big([\beta(x), y]+[x, \beta(y)]+\lambda[x, y]\Big),$
\item [3)] an {\bf element of centroid} \hspace{0,25cm}  if
$\beta([x, y])=[\beta(x), y]=[x, \beta(y)]$,
\item [4)] a {\bf Reynolds operator} \hspace{0,25cm} if
$\beta([x, y])=\beta\Big([\beta(x), y]+[x, \beta(y)]-[\beta(x), \beta(y)]\Big),$
\end{enumerate}
for any $x, y\in L$.
\end{definition}

\subsection{Ternary Leibniz algebras}
All the content of this subsection comes from \cite{IS}, where we have considered the linear vector space instead of graded vector space. 
\begin{definition}
 A ternary Leibniz algebra $L$ is said to be a ternary Leibniz  algebra if the bracket satisfies
 the following identity :
\begin{eqnarray}
 [[x, y, z], t, u]=[x, y, [z, t, u]]+[x, [y, t, u], z]+[[x, t, u], y, z]\label{lci}
\end{eqnarray}
for any $x, y, z, t, u\in L$.
\end{definition}

\begin{example}
 Let $A$ be a commutative associative algebra and $L$ a ternary Leibniz algebra. Then, for any $a, b, c\in A$, 
$x, y, z\in L$,  the bracket
\begin{eqnarray}
 \{a\otimes x, b\otimes y, c\otimes z\}=abc\otimes[x, y, z],
\end{eqnarray}
makes $A\otimes L$ into a ternary Leibniz algebra.
\end{example}

\begin{example}
 Let $L$ be a ternary Leibniz  algebra and put $L'=\mathbb{K}[t, t^{-1}]\otimes L$; $L'$ can be considered as a vector space over 
Laurent polynomials with coefficient in the ternary Leibniz algebra $L$.  Taking, for any $f(t), g(t), h(t)\in\mathbb{K}[t, t^{-1}]$ and 
$x, y, z\in L$,
\begin{eqnarray}
 [f(t)\otimes x, g(t)\otimes y, h(t)\otimes z]'=f(t)g(t)h(t)\otimes [x, y, z],\nonumber
\end{eqnarray}
we endow $L'$ with a structure of ternary Leibniz  algebra.
\end{example}

\begin{example}
 If $(L_1, [-, -, -]_1)$ and $(L_2,  [-, -, -]_2)$ are two  ternary Leibniz  algebras, the direct sum $L\oplus L'$ is also a ternary Leibniz  
algebra with respect to the operation :
\begin{eqnarray}
[x_1\oplus x_2, y_1\oplus y_2, z_1\oplus z_2]&:=&[x_1, y_1, z_1]_1\oplus[x_2, y_2, z_2]_2\nonumber,
\end{eqnarray}
for any $x_i, y_i, z_i\in L_i, i=1, 2$.
\end{example}

The following result asserts that one may associate a ternary Leibniz algebra to a Leibniz algebra. It will be useful later.
\begin{theorem}\label{ll3}
Let $(L, [-, -])$ be a Leibniz algebra. Then 
$$\bar L=(L, \{x, y, z\}:=[x, [y, z]]),$$
is a ternary Leibniz algebra, for  any $x, y, z\in L$.
\end{theorem}
\begin{proof}
 Applying twice relation (\ref{cpa}), for any $x, y, z\in L$, we have
\begin{eqnarray}
 &&\{\{x, y, z\}, t, u\}=\{[x, [y, z]], t, u\}=[[x, [y, z]], [t, u]]=\nonumber\\
&&\qquad=[x,  [[y, z], [t, u]]]+[[x, [t, u]]], [y, z]]\nonumber\\
&&\qquad=[x, [y, [z, [t, u]]]]+[x, [[y, [t, u]], z]]+[[x, [t, u]]], [y, z]]\nonumber\\
&&\qquad=\{x, y, [z, [t, u]]\}+[x, [\{y, t, u\}, z]]+[\{x, t, u\}, [y, z]]\nonumber\\
&&\qquad=\{x, y, \{z, t, u\}\}+\{x, \{y, t, u\}, z\}+\{\{x, t, u\}, y, z\}.\nonumber
\end{eqnarray}
This completes the proof.
\end{proof}

Here is the converse of Theorem \ref{ll3}.
\begin{proposition}\label{ltp}
 Let $(A, [-, -, -], \varepsilon)$ be a ternary Leibniz  algebra. Then $A\otimes A$ has a Leibniz algebra structure
 for the structure maps defined by
\begin{eqnarray}
\{x\otimes y, x'\otimes y'\}&:=&x\otimes[y, x', y']+[x, x', y']\otimes y.\nonumber
\end{eqnarray}
\end{proposition}

\begin{corollary}
 Let $(A, [-, -])$ be a Leibniz algebra. Then $A\otimes A$ is also a Leibniz  algebra
with respect to the operation
\begin{eqnarray}
\{x\otimes y, x'\otimes y'\}:=x\otimes[y, [x', y']]+[x, [x', y']]\otimes y.\nonumber
\end{eqnarray}
\end{corollary}

Now we recall modules on ternary Leibniz algebras.

\begin{definition}\label{mde}
 A bimodule $M$ over a ternary Leibniz  algebra $(A, \cdot, [-, -, -])$ is the given of three  trilinear applications
$$l_1 : M\otimes L\otimes L\rightarrow M, \quad l_2 : L\otimes M\otimes L\rightarrow M, \quad
l_3 : L\otimes L\otimes M\rightarrow M$$
satisfying the following sets of axioms
\begin{eqnarray}
 l_1(l_1(m, x, y), z, t)&=&l_1(m, x, [y, z, t])+l_1(m, [x, z, t], y)+l_1(l_1(m, z, t), x, y),\label{lpc3a1}\\
l_1(l_3(x, m, y), z, t)&=&l_2(x, m, [y, z, t])+l_2(x, l_1(m, z, t), y)+l_2([x, z, t], m, y),\\
l_1(l_3(x, y, m), z, t)&=&l_3(x, y, l_1(m, z, t))+l_3(x, [y, z, t], m)+l_3([x, z, t], y, m),\\
l_2([x, y, z], m, t)&=&l_3(x, y, l_2(z, m, t))+l_2(x, l_2(y, m, t), z)+l_1(l_2(x, m, t), y, z),\\
l_3([x, y, z], t, m)&=&l_3(x, y, l_3(z, t, m)+l_2(x, l_3(y, t, m), z)+l_1(l_3(x, t, m), y, z), \label{lpc3a5}
\end{eqnarray}
for all $x, y, z, t\in L, m\in M$.
\end{definition}
\begin{example}
 Any ternary Leibniz is a bimodule over it self.
\end{example}

The proof of the following proposition is straightforward.
\begin{proposition}
 Let $M$ and $N$ be two bimodules over a ternary Leibniz $L$.  Then, $M\oplus N$ is also a bimodule over $L$ with respect to 
componentwise operation.
\end{proposition}

The below proposition connects over Leibniz algebras to bimodules on ternary Leibniz algebras.
\begin{proposition}
 Let $(M, l, r)$ be a bimodule over a Leibniz  algebra $(L, [-, -])$. Let us define
\begin{eqnarray}
 \bar l_1(m, x, y):=r(m, [x,  y]),\quad \bar l_2(x, m, y):=l(x, r(m, y)),\quad \bar l_3(x, y, m):=l(x, l(y, m)),
\end{eqnarray}
for all $x, y\in L$ and $m\in M$.
Then $(M, \bar l_1, \bar l_2, \bar l_3)$ is a bimodule over the ternary Leibniz algebra  $\bar L$ (as in Theorem \ref{ll3}).
\end{proposition}

\section{Operators on compatible ternary Leibniz algebras}
In this section, we study the invariance of some operators on the (compatible) ternary Leibniz algebra induced by (compatible)
 Leibniz algebra, makes the relationship among these operators and  establish the ternary Leibniz algebra strucutre induced by either
 averaging operators or $O$-operator.

\begin{definition}
 A compatible ternary Leibniz algebra is a triple $(L, [-, -, -]_1, [-, -, -]_2)$ where  $(L, [-, -, -]_1)$
and $(L, [-, -, -]_2)$ are two ternary Leibniz algebras such that :
\begin{eqnarray}
 [[x, y, z]_1, t, u]_2+[[x, y, z]_2, t, u]_1&=&[x, y, [z, t, u]_1]_2+[x, y, [z, t, u]_2]_1+[x, [y, t, u]_1, z]_2\nonumber\\
&&+[x, [y, t, u]_2, z]_1+[[x, t, u]_1, y, z]_2+[[x, t, u]_2, y, z]_1, \label{cc}
\end{eqnarray}
for any $x, y, z\in L$.
\end{definition}

The proof of the next proposition is straightforward.
\begin{proposition}
 The triple $(L, [-, -, -]_1, [-, -, -]_2)$ is a compatible ternary Leibniz algebra if and only if the new bilinear map 
$$[x, y, z]=k_1[-, -, -]_1+k_2[-, -, -]_2$$
define a ternary Leibniz structure on $L$.
\end{proposition}
\subsection{Invariance of operators on compatible ternary Leibniz algebras}
Along this subsection, an operator on two compatible (ternary) Leibniz algebras is an operator for both the two structures.
\begin{definition}
\begin{enumerate}
 \item [1)]
Let $(L, [-, -])$ be a Leibniz algebra. A
linear map $D\in End(L)$ is said to be a central derivation if 
$$D([x, y])=[D(x), y]=[x, D(y)]=0,$$
for all $x, y\in L$.
\item [2)]
Let $(L, [-, -, -])$ be a ternary Leibniz algebra. A
linear map $D\in End(L)$ is said to be a central derivation if
$$D([x, y, z]) =[D(x), y, z]=[x, D(y), z]=[x, y, D(z)]=0,$$
for all $x, y\in L$.
\end{enumerate}
\end{definition}

\begin{lemma}
 Let $D : L\rightarrow L$ be a central derivation of a Leibniz algebra $(L, [-, -])$.
Then $D$ is also a central derivation of the associated ternary Leibniz algebra. 
\end{lemma}
\begin{proof}
 For any $x, y, z\in L$, 
\begin{eqnarray}
D(\{x, y, z\})&=&D([x, [y, z]])=[D(x), [y, z]]=[x, D([y, z])]=[x, [D(y), z]]=[x, [y, D(z)]]=0\nonumber\\
&=&\{D(x), y, z\}=\{x, D(y), z\}=\{x, y, D(z)\}=0.\nonumber
\end{eqnarray}
Which implies $D$ is a central derivation of the ternary Leibniz algebra.
\end{proof}

\begin{proposition}
  Let $(L, [-, -]_1)$ and $(L, [-, -]_2)$ be two compatible Leibniz algebras such that the associated ternary Leibniz algebras be compatibles.
Let $D : L\rightarrow L$ be a central derivation of the compatible Leibniz algebra $(L, [-, -]_1, [-, -]_2)$ i.e., for any $x, y\in L$,
$$D([x, y]_i)=[D(x), y]_i=[x, D(y)]_i=0, i=1, 2.$$
Then $D$ is also a central derivation of the compatible ternary Leibniz algebra $\bar L$. 
\end{proposition}
\begin{definition}
\begin{enumerate}
 \item [1)]
Let $(L, [-, -])$ be a Leibniz algebra. A
linear map $D\in End(L)$ is said to be a generalized derivation if there exists two linear maps $D_1, D_2\in End(L)$  such that
$$D''([x, y]) =[D(x), y]+[x, D'(y)],$$
for all $x, y\in L$.
\item [2)]
Let $(L, [-, -, -])$ be a ternary Leibniz algebra. A
linear map $D\in End(L)$ is said to be a generalized derivation if there
exists three linear maps $D_1, D_2, D_3\in End(L)$  such that
$$D'''([x, y, z]) =[D(x), y, z]+[x, D'(y), z]+[x, y, D''(z)],$$
for all $x, y\in L$.
\end{enumerate}
\end{definition}
\begin{lemma}
 Let $D : L\rightarrow L$ be a generalized derivation on a Leibniz algebra $(A, [-, -])$ such that
 \begin{eqnarray}
  D''([x, y])&=&[D(x), y]+[x, D'(y)],\nonumber\\
D'''([x, y])&=&[D(x), y]+[x, D''(y)],\nonumber
 \end{eqnarray}
for any $x, y\in L$. Then, $D$ is a generalized derivation of the ternary Leibniz algebra $\bar L$.  
\end{lemma}
\begin{proof}
 For any $x, y, z\in L$, 
\begin{eqnarray}
 D'''(\{x, y, z\})&=&D'''([x, [y, z]])=[D(x), [y, z]]+[x, D''([y, z])]\nonumber\\
&=&[D(x), [y, z]]+ [x, [D(y), z]+ [y, D'(z)]]\nonumber\\
&=&[D(x), [y, z]]+ [x, [D(y), z]]+[x, [y, D'(z)]]\nonumber\\
&=&\{D(x), y, z\}+ \{x, D(y), z\}+\{x, y, D'(z)\}\nonumber.
\end{eqnarray}
Therefore, $D$ is a generalized derivation of $\bar L$.
\end{proof}

\begin{proposition}
  Let $(L, [-, -]_1)$ and $(L, [-, -]_2)$ be two compatible Leibniz algebras such that the associated ternary Leibniz algebras be compatible.
 Let $D : L\rightarrow L$ be a generalized derivation on the compatible Leibniz algebra $(A, [-, -]_1, [-, -]_2)$ such that
 \begin{eqnarray}
  D''([x, y]_i)&=&[D(x), y]_i+[x, D'(y)]_i,\nonumber\\
D'''([x, y]_i)&=&[D(x), y]_i+[x, D''(y)]_i,\nonumber
 \end{eqnarray}
for any $x, y\in L$. Then, $D$ is a generalized derivation of the compatible ternary Leibniz algebra $\bar L$.  
\end{proposition}
\begin{definition}
\begin{enumerate}
 \item [1)]
Let $(L, [-, -])$ be a Leibniz algebra. A
linear map $D\in End(L)$ is said to be a quasiderivation if there exists a linear map $D'\in End(L)$  such that
$$D'([x, y]) =[D(x), y]+[x, D(y)],$$
for all $x, y\in L$.
\item [2)]
Let $(L, [-, -, -])$ be a ternary Leibniz algebra. A
linear map $D\in End(L)$ is said to be a quasiderivation if there
exists a linear map $D'\in End(L)$  such that
$$D'([x, y, z]) =[D(x), y, z]+[x, D(y), z]+[x, y, D(z)],$$
for all $x, y\in L$.
\end{enumerate}
\end{definition}

\begin{lemma}\label{qsd}
 Let $D, D' : L\rightarrow L$ be two quasi-derivations on a Leibniz algebra $(A, [-, -])$ such that
 \begin{eqnarray}
  D'([x, y])&=&[D(x), y]+[x, D(y)],\nonumber\\
D''([x, y])&=&[D'(x), y]+[x, D'(y)],\nonumber
 \end{eqnarray}
for any $x, y\in L$. Then, $D$ is a quasi-derivation of the ternary Leibniz algebra.  
\end{lemma}
\begin{proof}
 For all $x, y\in A$, 
 \begin{eqnarray}
D''(\{x, y, z\})&=&D''([[x, y], z])=[D'[x, y], z]+[[x, y], D'(z)]\nonumber\\
&=&[[D(x), y]+ [x, D(y)], z]+[[x, y], D'(z)]\nonumber\\
&=&\{D(x), y, z\}+ \{x, D(y), z\}+\{x, y, D'(z)\}\nonumber.
\end{eqnarray}
This proves the assertion.
\end{proof}

\begin{proposition}
  Let $(L, [-, -]_1)$ and $(L, [-, -]_2)$ be two compatible Leibniz algebras such that the associated ternary Leibniz algebras be compatible.
 Let $D, D' : L\rightarrow L$ be two quasi-derivations on the compatible Leibniz algebra $(A, [-, -]_1, [-, -]_2)$ such that
 \begin{eqnarray}
  D'([x, y]_i)&=&[D(x), y]_i+[x, D(y)]_i,\nonumber\\
D''([x, y]_i)&=&[D'(x), y]_i+[x, D'(y)]_i,\nonumber
 \end{eqnarray}
for any $x, y\in L$. Then, $D$ is a quasi-derivation of the compatible ternary Leibniz algebra $\bar L$.  
\end{proposition}
\begin{proof}
 The proof follows from Lemma \ref{qsd} by direct calculation.
\end{proof}

\begin{definition}
 A linear map $R : L\rightarrow L$ on a ternary Leibniz algebra is called 
\begin{enumerate}
 \item [1)]
a {\bf Rota-Baxter operator of weight $\lambda\in\mathbb{K}$}, if 
\begin{eqnarray}
 &&[R(x), R(y), R(z)]=R\Big([R(x), R(y), z]+[R(x), y, R(z)]+[x, R(y), R(z)]\nonumber\\
&&\qquad\qquad\qquad\qquad+\lambda[R(x), y, z]+\lambda[x, R(y), z]+\lambda[x, y, R(z)]+\lambda^2
[x, y, z]\Big),\label{rotid}
\end{eqnarray}
\item [2)]
 a {\bf derivation of weight $\lambda\in\mathbb{K}$} if 
\begin{eqnarray}
 &&R([x, y, z])=[R(x), R(y), z]+[R(x), y, R(z)]+[x, R(y), R(z)]\nonumber\\
&&\qquad\qquad\qquad\qquad+\lambda[R(x), y, z]+\lambda[x, R(y), z]+\lambda[x, y, R(z)]+\lambda^2[x, y, z],\label{rot}
\end{eqnarray}
\end{enumerate}
for any $x, y, z\in L$,
\end{definition}
\begin{proposition}
 Let $(L, [-, -, -])$ be a ternary Leibniz algebra and $R: L\rightarrow L$ an inversible Rota-Baxter operator of weight $\lambda\in\mathbb{K}$ on 
$L$. Then, $R^{-1}: L\rightarrow L$ is a derivation of weight $\lambda$.  
\end{proposition}
\begin{proof}
 It follows from Rota-Baxter identity (\ref{rotid}) by putting $R(x)=x', R(y)=y', R(z)=z'$.
\end{proof}
 We continue to define operators on ternary Leibniz algebras.
\begin{definition}\label{nod}
 A linear map $N : L\rightarrow L$ on a ternary Leibniz algebra is called a {\bf Nijenhuis operator} if
\begin{eqnarray}
 &&[N(x), N(y), N(z)]=N\Big([N(x), N(y), z]+[N(x), y, N(z)]+[x, N(y), N(z)]\Big)\nonumber\\
&&\qquad\qquad\qquad\qquad-N^2\Big([N(x), y, z]+[x, N(y), z]+[x, y, N(z)]\Big)+N^3([x, y, z])\Big),\label{rot}
\end{eqnarray}
for any $x, y, z\in L$.
\end{definition}

The following propositions connect Rota-Baxter operators to derivations (see Proposition 3.1 \cite{XDC} for pre-Lie algebras) 
and Nijenhuis operators (see Proposition 2.11, \cite{BR} for Leibniz algebras). 
\begin{proposition}
 Let $(L, [-, -, -])$ be a ternary Leibniz algebra and $N : L\rightarrow L$ be a linear operator.
\begin{enumerate}
 \item [1)] If $N^2=0$, then $N$ is a Nijenhuis operator if and only if $N$ is a Rota-Baxter operator of weight $0$.
\item [2)] If $N^2=N$ then $N$ is a Nijenhuis operator if $N$ is a Rota-Baxter operator of weight $-1$.
\end{enumerate}
\end{proposition}
\begin{proof}
\begin{enumerate}
 \item [1)]
As $N^2=0$ and $N$ is a Nijenhuis operator, we have 
\begin{eqnarray}
 [N(x), N(y), N(z)]&=&N\Big([N(x), N(y), z]+[N(x), y, N(z)]+[x, N(y), N(z)]\Big)\nonumber\\
 &&\quad-N^2\Big([N(x), y, z]+[x, N(y), z]+[x, y, N(z)]\Big)+N^3([x, y, z]),\nonumber\\
&=&N\Big([N(x), N(y), z]+[N(x), y, N(z)]+[x, N(y), N(z)]\Big).\nonumber
\end{eqnarray}
Which means that $N$ is a Rota-Baxter operator of weight $0$. The converse is proved by going up the egalities.
\item [2)] Since $N^2=N$ and $N$ is a Nijenhuis operator, we have 
\begin{eqnarray}
 [N(x), N(y), N(z)]&=&N\Big([N(x), N(y), z]+[N(x), y, N(z)]+[x, N(y), N(z)]\Big)\nonumber\\
 &&\quad-N^2\Big([N(x), y, z]+[x, N(y), z]+[x, y, N(z)]\Big)+N^3([x, y, z]),\nonumber\\
 &=&N\Big([N(x), N(y), z]+[N(x), y, N(z)]+[x, N(y), N(z)]\Big)\nonumber\\
 &&\quad-N\Big([N(x), y, z]+[x, N(y), z]+[x, y, N(z)]\Big)+N([x, y, z]),\nonumber\\
&=&N\Big([N(x), N(y), z]+[N(x), y, N(z)]+[x, N(y), N(z)]\nonumber\\
 &&\quad-[N(x), y, z]+[x, N(y), z]+[x, y, N(z)]+[x, y, z]\Big).\nonumber
\end{eqnarray}
Which means that $N$ is a Rota-Baxter operator of weight $0$. The converse is proved by going up the egalities.
\end{enumerate}
\end{proof}
\begin{definition}
 A linear map $\theta : L\rightarrow L$ is said to be an element of centroid of $L$ if, for any $x, y, z\in\mathcal{H}(L)$,
\begin{eqnarray}
 \theta([x, y, z])=[\theta(x), y, z]=[x, \theta(y), z]=[x, y, \theta(z)].\label{ece}
\end{eqnarray}
\end{definition}
\begin{remark}
The sum of two element of centroids is also another one.
\end{remark}
\begin{example}
For any $\lambda\in\mathbb{K}$ and any $x\in L$, the linear map $\theta(x)=\lambda x$ is an element of centroid of $L$.
\end{example}
\begin{example}
 Let $L$ be a ternary Leibniz color algebra and $\theta : L\rightarrow L$ an element of centroid of $L$. Then $L$ becomes a ternary Leibniz
 color algebra in each of the following brackets :
 $$\{x, y, z\}_1=[\theta(x), y, z], \quad\{x, y, z\}_2=[\theta(x), \theta(y), z],\quad\{x, y, z\}_3=[\theta(x), \theta(y), \theta(z)],$$
 for any $x, y, z\in L$.
\end{example}
\begin{definition}
 A linear map $R : L\rightarrow L$ on a ternary Leibniz algebra is called a {\bf Reynolds operator} if
\begin{eqnarray}
 [R(x), R(y), R(z)]\!=\!R\Big([R(x), R(y), z]\!+\![R(x), y, R(z)]\!+\![x, R(y), R(z)]\!-\![R(x), R(y), R(z)]\Big), \label{ret}
\end{eqnarray}
for any $x, y, z\in L$.
\end{definition}

The proof of the following proposition is straightforward.
\begin{theorem}\label{ops}
 Let $\varphi$ ($D, N, \beta, R$, $P$ respectively) be an element of centroid (a derivation, a Nijenhuis operator, an averaging operator, 
a Rota-Baxter operator, a Reynolds operator respectively) on a Leibniz algebra $L$. Then, $\varphi$ ($D, N, \beta, R$, $P$ respectively)
 is also  an element of centroid (a derivation, a Nijenhuis operator, an averaging operator, a Rota-Baxter operator, a Reynolds operator respectively) on the associated ternary 
Leibniz algebra.
\end{theorem}
\begin{proof}
 We only prove for element of centroid. So, for any $x, y, z\in L$,
\begin{eqnarray}
 \theta([x, y, z])=\theta([x, [y, z]])=[\theta(x), [y, z]]=[x, \theta([y, z])]=[x, [\theta(y), z]]=[x, [y, \theta(z)]]\nonumber.
\end{eqnarray}
Thus, $$\theta([x, y, z])=[\theta(x), y, z]=[x, \theta(y), z]=[x, y, \theta(z)].$$
The other are done in the same way.
\end{proof}
\subsection{Averaging operators on compatible ternary Leibniz algebras}
\begin{definition}
 A  linear map $\beta : L\rightarrow L$ is said to be an averaging operator on $L$ if, for any $x, y, z\in L$,
\begin{eqnarray}
 \beta([\beta(x), \beta(y), z]=\beta([\beta(x), y, \beta(z)])=\beta([x, \beta(y), \beta(z)])=[\beta(x), \beta(y), \beta(z)].\label{ece}
\end{eqnarray}
\end{definition}

\begin{lemma}\label{avt}
 Let $(L, [-, -, -, \beta)$ be an injective averaging ternary Leibniz algebra. Then, $L$, endowed with the new multiplication 
$[-, -, -] : L\times L\times L\rightarrow L$ defined by 
\begin{eqnarray}
 [x, y, z]_\beta:=[\beta(x), \beta(y), z],
\end{eqnarray}
for all $x, y, z\in L$, makes $L$ into a ternary Leibniz algebra.
\end{lemma}
\begin{proof}
 For any $x, y, z\in L$,
\begin{eqnarray}
&&\qquad\beta\Big([[x, y, z]_\beta, t, u]_\beta-[x, y, [z, t, u]_\beta]_\beta-[x, [y, t, u]_\beta, z]_\beta-[[x, t, u]_\beta, y, z]_\beta\Big)\nonumber\\
&&=\beta\Big([\beta([\beta(x), \beta(y), z]), \beta(t), u]-[\beta(x), \beta(y), [\beta(z), \beta(t), u]]-
[\beta(x), \beta([y, t, u]), z]-[\beta([x, t, u]), \beta(y), z]\Big)\nonumber\\
&&=[\beta([\beta(x), \beta(y), z]), \beta(t), \beta(u)]-[\beta(x), \beta(y), \beta([\beta(z), \beta(t), u])]\nonumber\\
&&\qquad-[\beta(x), \beta([\beta(y), \beta(t), u]), \beta(z)]-[\beta([\beta(x), \beta(t), u]), \beta(y), \beta(z)]\nonumber\\
&&=[[\beta(x), \beta(y), \beta(z)], \beta(t), \beta(u)]-[\beta(x), \beta(y), [\beta(z), \beta(t), \beta(u)]]\nonumber\\
&&\qquad-[\beta(x), [\beta(y), \beta(t), \beta(u)]), \beta(z)]-[[\beta(x), \beta(t), \beta(u)]), \beta(y), \beta(z)]\nonumber.
\end{eqnarray}
The right hand side vanishes by ternary Leibniz rule, and the conclusion follows from injectivity.
\end{proof}

\begin{proposition}
Let $(L, [-, -, -]_1, \beta_1)$ and $(L, [-, -, -]_2, \beta_2)$ be two commuting
injective averaging ternary Leibniz algebras  such that
$$[\beta_1(x), \beta_1(y), \beta_1(z)]_1=[\beta_2(x), \beta_2(y), \beta_2(z)]_2,$$
for all $x, y, z\in L$. Then, with notation of Lemma \ref{avt}, $(L, [-, -, -]_{\beta_1}, [-, -, -]_{\beta_2})$ is a compatible ternary Leibniz algebra.
\end{proposition}
\begin{proof}
 It follows from direct computation of 
\begin{eqnarray}
 &&\beta_1\beta_2\Big([[x, y, z]_{\beta_1}, t, u]_{\beta_2}+[[x, y, z]_{\beta_2}, t, u]_{\beta_1}
-[x, y, [z, t, u]_{\beta_1}]_{\beta_2}-[x, y, [z, t, u]_2]_{\beta_1}-[x, [y, t, u]_{\beta_1}, z]_{\beta_2}\nonumber\\
&&\hspace{5cm}-[x, [y, t, u]_{\beta_2}, z]_{\beta_1}-[[x, t, u]_{\beta_1}, y, z]_{\beta_2}-[[x, t, u]_{\beta_2}, y, z]_{\beta_1}\Big),\nonumber
\end{eqnarray}
by using the assumption.
\end{proof}

\begin{remark}
 In the above proposition, we start with two injective averaging ternary Leibniz algebras which are not compatible at the begining.
\end{remark}

\begin{proposition}
 Let $(A, [-, -, -]_1, [-, -, -]_2)$ be a compatible ternary Leibniz algebra and $\beta : A\rightarrow A$ be an injective averaging operator
on $A$ for both $[-, -, -]_1$ and $[-, -, -]_2)$. Then, $A^\beta :=(A, [-, -, -]_1^\beta, [-, -, -]_2^\beta)$ is a compatible ternary Leibniz
algebra, where $[-, -, -]_i^\beta=[\beta(-), \beta(-), -]_i, i=1,2$. Moreover, $\beta : (A, [-, -, -]^\beta)\rightarrow (A, [-, -, -])$ is a
morphism of ternary Leibniz algebra, with $[-, -, -]=[-, -, -]_1+[-, -, -]_2$.
\end{proposition}
\begin{proof}
It is clear, from Lemma \ref{avt}, that $A$ is a ternary Leibniz algebra for each of the bracket $[-, -, -]_1^\beta$ and $[-, -, -]_2^\beta$.
Then, for any $x, y, z\in A$,
\begin{eqnarray}
 &&\qquad\beta\Big([[x, y, z]_1^\beta, t, u]_2^\beta+[[x, y, z]_2^\beta, t, u]_1^\beta
-[x, y, [z, t, u]_1^\beta]_2^\beta- [x, y, [z, t, u]_2^\beta]_1^\beta-[x, [y, t, u]_1^\beta, z]_2^\beta\nonumber\\
&&\qquad-[x, [y, t, u]_2^\beta, z]_1^\beta-[[x, t, u]_1^\beta, y, z]_2^\beta-[[x, t, u]_2\beta, y, z]_1^\beta\Big)=\nonumber\\
&&=\beta\Big([\beta([\beta(x), \beta(y), z]_1), \beta(t), u]_2+[\beta([\beta(x), \beta(y), z]_2), \beta(t), u]_1
-[\beta(x), \beta(y), [\beta(z), \beta(t), u]_1]_2\nonumber\\
&&\qquad- [\beta(x), \beta(y), [\beta(z), \beta(t), u]_2]_1
-[\beta(x), \beta([\beta(y), \beta(t), u]_1), z]_2-[\beta(x), \beta([\beta(y), \beta(t), u]_2), z]_1\nonumber\\
&&\qquad-[\beta([\beta(x), \beta(t), u]_1, \beta(y), z]_2-[\beta([\beta(x), \beta(t), u]_2, \beta(y), z]_1\Big)\nonumber\\
&&=[\beta([\beta(x), \beta(y), z]_1), \beta(t), \beta(u)]_2+[\beta([\beta(x), \beta(y), z]_2), \beta(t), \beta(u)]_1
-[\beta(x), \beta(y), \beta([\beta(z), \beta(t), u]_1)]_2\nonumber\\
&&\qquad- [\beta(x), \beta(y), \beta([\beta(z), \beta(t), u]_2)]_1
-[\beta(x), \beta([\beta(y), \beta(t), u]_1, \beta(z)]_2-[\beta(x), \beta([\beta(y), \beta(t), u]_2), \beta(z)]_1\nonumber\\
&&\qquad-[\beta([\beta(x), \beta(t), u]_1), \beta(y), \beta(z)]_2-[\beta([\beta(x), \beta(t), u]_2, \beta(y), \beta(z)]_1\nonumber\\
&&=[[\beta(x), \beta(y), \beta(z)]_1, \beta(t), \beta(u)]_2+[[\beta(x), \beta(y), \beta(z)]_2), \beta(t), \beta(u)]_1
-[\beta(x), \beta(y), [\beta(z), \beta(t), \beta(u)]_1)]_2\nonumber\\
&&\qquad- [\beta(x), \beta(y), [\beta(z), \beta(t), \beta(u)]_2]_1
-[\beta(x), [\beta(y), \beta(t), \beta(u)]_1, \beta(z)]_2-[\beta(x), [\beta(y), \beta(t), \beta(u)]_2), \beta(z)]_1\nonumber\\
&&\qquad-[[\beta(x), \beta(t), \beta(u)]_1, \beta(y), \beta(z)]_2-[[\beta(x), \beta(t), \beta(u)]_2, \beta(y), \beta(z)]_1\nonumber
\end{eqnarray}
The right hand side vanishes by compatibility condition, and the conclusion follows from injectivity.\\
Now,
\begin{eqnarray}
\beta([x, y, z]^\beta)&=&\beta\Big([x, y, z]_1^\beta+[x, y, z]_2^\beta\Big)
=\beta\Big([\beta(x), \beta(y), z]_1+[\beta(x), \beta(y), z]_2\Big)\nonumber\\
&=&[\beta(x), \beta(y), \beta(z)]_1+[\beta(x), \beta(y), \beta(z)]_2\nonumber\\
&=&[\beta(x), \beta(y), \beta(z)]\nonumber.
\end{eqnarray}
This achieves the proof.
\end{proof}

\begin{corollary}
 Let $(L, [-, -]_1, [-, -]_2)$ be a compatible Leibniz algebra such that the associated ternary Leibniz algebra is compatible and 
$\beta : L\rightarrow L$ an injective averaging operator on $L$. Then, 
$(L, [-, [-, -]_1^\beta]_1^\beta, [-, [-, -]_2^\beta]_2^\beta)$ is a compatible ternary Leibniz algebra.
\end{corollary}
\begin{proof}
 It follows from that $[-, -]^\beta=[\beta(-), -]$ is a Leibniz algebra and
\begin{eqnarray}
 [x, y, z]_i^\beta=[x, [y, z]_i^\beta]_i^\beta=[\beta(x), [\beta(y), z]_i]_i=[\beta(x), \beta(y), z]_i,\nonumber
\end{eqnarray}
for any $x, y, z\in L$.
\end{proof}

\begin{proposition}
 Let us denote $\mathcal C_c(L)$ the category of all compatible ternary Leibniz structures on the vector space $L$ 
(the category $\mathcal C_c(L)$ contains only one object, $L$). On $\mathcal C_c(L)$, let us define the relation 
$$(L, [-, -]_1)\mathcal R (L, [-, -]_2)\quad \mbox{if and only if}\quad (L, [-, -]_1, [-, -]_2)\quad \mbox{is a ternary Leibniz algebra}$$
The relation $\mathcal R$ is an equivalence relation on $\mathcal C_c(L)$.
\end{proposition}
\begin{proof}
 It is trivial by construction of $\mathcal C_c(L)$.
\end{proof}

\subsection{Rota-Baxter modules on compatible ternary Leibniz algebras}
The following definition is inspired from the work on bimodules and Rota-Baxter relation (\cite{BB}). 
\begin{definition}
 Let $(L, [-, -, -])$ be a ternary algebra, $R : L\rightarrow L$ a Rota-Baxter operator on $L$ and $M$ a bimodule on $L$ with structure maps
$l_1, l_2, l_3$ and $R_M : M\rightarrow M$ a linear map on $M$. We say that $(M, R_M)$ is a Rota-Baxter bimodule on $L$ if, for any $x, y\in L$,
$m\in M$, we have 
\begin{eqnarray}
 l_1(R_M(m), R(x), R(y))&=&R_M\Big(l_1(R_M(m), R(x), y)+l_1(m, R(x), R(y))+l_1(R_M(m), x, R(y))\Big),\\
l_2(R(x), R_M(m), R(y))&=&R_M\Big(l_2(R(x),R_M(m), y)+l_2(x, R_M(m), R(y))+l_2(R(x), m, R(y))\Big),\\
l_3(R(x), R(y), R_M(m))&=&R_M\Big(l_3(R(x), y, R_M(m))+l_3(x, R(y), R_M(m))+l_3(R(x), R(y), m)\Big).
\end{eqnarray}
\end{definition}

\begin{theorem}
 Let $(M, l_1, l_2, l_3, R_M)$ be a Rota-Baxter bimodule over the Rota-Baxter ternary Leibniz algebra $(L, [-, -, -], R)$.
Let us define the three linear maps $l_1', l_2', l_3'$ as follows
 \begin{eqnarray}
 l_1'(m, x, y)&=&l_1(R_M(m), R(x), y)+l_1(m, R(x), R(y))+l_1(R_M(m), x, R(y)),\nonumber\\\
l_2'(x, m, y)&=&l_2(R(x),R_M(m), y)+l_2(x, R_M(m), R(y))+l_2(R(x), m, R(y)),\nonumber\\
l_3'(x, y, m)&=&l_3(R(x), y, R_M(m))+l_3(x, R(y), R_M(m))+l_3(R(x), R(y), m).\nonumber
\end{eqnarray}
Then $(M, l'_1, l'_2, l_3')$ is a bimodule over the ternary Leibniz algebra $L_R=(L, [-, -, -]_R)$.
\end{theorem}
 \begin{proof}
  We only prove the first axiom. The other being done similarly. Then, for any $x, y, z, t\in L$ and any $m\in M$, we have
\begin{eqnarray}
 &&\qquad l'_1(l'_1(m, x, y), z, t)-l'_1(m, x, [y, z, t]_R)-l'_1(m, [x, z, t]_R, y)-l'_1(l'_1(m, z, t), x, y)=\nonumber\\
&&=l_1(R_Ml'_1(m, x, y), R(z), t)+l_1(l_1'(m, x, y), R(z), R(t))+l_1(R_Ml'_1(m, x, y), z, R(t))\\
&&\quad-l_1(R_M(m), R(x), [y, z, t]_R)-l_1(m, R(x), R([y, z, t]_R))-l_1(R_M(m), x, R([y, z, t]_R))\nonumber\\
&&\quad-l_1(R_M(m), R([x, z, t]_R), y)-l_1(m, R([x, z, t]_R), R(y))-l_1(R_M(m), [x, z, t]_R, R(y))\nonumber\\
&&\quad-l_1(R_M(l_1'(m, z, t)), R(x), y)-l_1(l_1'(m, z, t), R(x), R(y))-l_1(R_Ml_1'(m, z, t), x, R(y))\nonumber\\
&&=l_1(l_1(R_M(m), R(x), R(y)), R(z), t)+l_1\Big(l_1(R_M(m), R(x), y)+l_1(m, R(x), R(y))+l_1(R_M(m), x, R(y)), R(z), R(t)\Big)\nonumber\\
&&\qquad+l_1(l_1(R_M(m), R(x), R(y)), z, t)-l_1\Big(R_M(m), R(x), [R(y), R(z), t]+[R(y), z, R(t)]+[y, R(z), R(t)]\Big)\nonumber\\
&&\quad-l_1(m, R(x), [R(y), R(z), R(t)])-l_1(R_M(m), x, [R(y), R(z), R(t)])-l_1(R_M(m), [R(x), R(z), R(t)], y)\nonumber\\
&&\qquad-l_1(m, [R(x), R(z), R(t)], R(y))-l_1\Big(R_M(m), [R(x), R(z), t]+[R(x), z, R(t)]+[x, R(z), R(t)], R(y)\Big)\nonumber\\
&&\quad-l_1(l_1(R_M(m), R(z), t), y)-l_1\Big(l_1(R_M(m), R(z), t)+l_1(R_M(m), z, R(t))+l_1(m, R(z), R(t)), R(x), R(y)\Big)\nonumber\\
&&\qquad-l_1(l_1(R_M(m), R(z), R(t)), x, R(y))\nonumber.
\end{eqnarray}
The left hand side vanishes by axiom (\ref{lpc3a1}), and the conclusion follows.
 \end{proof}

\subsection{$O$-operators on compatible ternary Leibniz algebras}
\begin{definition}(\cite{IS})\label{dr}
 Let $L$ be a ternary Leibniz  algebra and $M$ a  vector space. A representation of $L$ on $M$ is the given of three  linear
maps $\lambda : L\times L\rightarrow End(M), \mu : L\otimes L\rightarrow End(M), $ and $\rho : L\otimes L\rightarrow End(M)$ such that :
\begin{eqnarray}
 \lambda({[x, y, z], t})(m)&=& \lambda({x, y})\lambda({z, t})(m)+\mu({x, z})\lambda({y, t})(m)+\rho({y, z})\lambda({x, t})(m),\\
\mu({[x, y, z], t})(m)&=&\lambda({x, y})\mu({z, t})(m)+\mu({x, z})\mu({y, t})(m)+\rho({y, z})\mu({x, t})(m),\\
\rho({z, t})\lambda({x, y})(m)&=&\lambda({x, y})\rho({z, t})(m)+\lambda({x, [y, z, t]})(m)+\lambda({[x, z, t], y})(m),\\
\rho({z, t})\mu({x, y})(m)&=&\mu({x, [y, z, t]})(m)+\mu({x, y})\rho({z, t})(m)+\mu({[x, z, t], y})(m),\\
\rho({z, t})\rho({x, y})(m)&=&\rho({x, [y, z, t]})(m)+\rho({[x, z, t], y})(m)+\rho({x, y})\rho({z, t})(m),
\end{eqnarray}
for any $x, y, z, t\in L, m\in M$.
\end{definition}

\begin{remark}
It is well kwnown that every  bimodule $M$ gives rise to a representation $(\lambda, \mu, \rho)$ of $L$ on $M$ via 
$\lambda({x, y})(m)=l_3(x, y, m)$, $\mu({x, y})(m)=l_2(x, m, y)$ and $\rho({x, y})(m)=l_1(m, x, y)$.
 Conversely, every representation $(\lambda, \mu, \rho)$ of $L$ on $M$, we define an $L$ bimodule structure on $M$ via
$l_3(x, y, m):=\lambda({x, y})(m)$,
$l_2(x, m, y):=\mu({x, y})(m)$ and $l_1(m, x, y):=\rho({x, y})(m)$. 
\end{remark}

\begin{definition}
Let $(L, [-, -, -])$ be a ternary Leibniz algebra, $(V, \lambda, \mu, \rho)$ a representation of $L$ on $V$ and $T : V\rightarrow L$ a linear 
map. We say that $T$ is an $O$-operator on $L$ with respect to the representation $(\lambda, \mu, \rho)$ if 
 \begin{eqnarray}
 [T(u), T(v), T(w)]=T\Big(\lambda(T(u), T(v))w+\mu(T(u), T(w))v+\rho(T(v), T(w))u\Big),
\end{eqnarray}
for any $u, v, w \in V$.
\end{definition}

\begin{definition}
Let $T, T'$ be two $O$-operators  on a ternary Leibniz algebra $(L, [-, -, -])$ with respect to a representation $(V, \lambda, \mu, \rho)$
 A morphism from $T'$ to $T$ consists of a ternary Leibniz algebra morphism $f : L\rightarrow L$ and a linear map $\varphi : V\rightarrow V$ such
that 
\begin{eqnarray}
 T\circ\varphi&=&f\circ T',\\
(\varphi\circ \lambda(x, y))(u)&=&\lambda(f(x), f(y))\varphi(u),\\
(\varphi\circ \mu(x, y))(u)&=&\mu(f(x), f(y))\varphi(u),\\
(\varphi\circ \rho(x, y))(u)&=&\rho(f(x), f(y))\varphi(u),
\end{eqnarray}
for any $x, y\in L, u\in V$.
\end{definition}

\begin{theorem}
Let $T: V\rightarrow L$ be an $O$-operator on $L$ with respect to the representation $(V, \lambda, \rho, \mu)$. Let us define the bracket
$[-, -, -]_T : V\times V\times V\rightarrow V$ as follows
\begin{eqnarray}
[u_1, u_2, u_3]_T= \lambda(T(u_1), T(u_2))(u_3)+\mu(T(u_1), T(u_3))(u_2)+\rho(T(u_2), T(u_3))(u_1)
\end{eqnarray}
for any $u_1, u_2, u_3\in V$. Then $(V, [-, -, -]_T)$ is a ternary Leibniz algebra.\\
 Let $T'$ be another $O$-operator on a ternary Leibniz algebra $(L, [-, -, -])$ with respect to a representation $(V, \lambda, \mu, \rho)$.
and $(f, \phi)$ a homomorphism from $T'$ to $T$. Then, $\phi$ is a homomorphism from ternary Leibniz algebra $(V, [-, -, -]_{T'})$ to
$(V, [-, -, -]_{T})$.
\end{theorem}
\begin{proof}
For any $u_i\in V, i=1,2,3,4,5$, we have
 \begin{eqnarray}
[[u_1, u_2, u_3]_T, u_4, u_5]_T&=&\lambda(T([u_1, u_2, u_3]_T), u_4)(u_5)+\mu(T([u_1, u_2, u_3]_T), T(u_5))(u_4)
+\rho(T(u_4), T(u_5))([u_1, u_2, u_3]_T)\nonumber\\
&=&\lambda(T(u_1), T(u_2), T(u_3), u_4)(u_5)+\mu(T(u_1), T(u_2), T(u_3), T(u_5))(u_4)\nonumber\\
&&+\rho(T(u_4), T(u_5))\Big(\lambda(T(u_1), T(u_2))(u_3)+\mu(T(u_1), T(u_3))(u_2)+\rho(T(u_2), T(u_3))(u_1)\Big)\nonumber\\
&=&\lambda(T(u_1), T(u_2), T(u_3), u_4)(u_5)+\mu(T(u_1), T(u_2), T(u_3), T(u_5))(u_4)\nonumber\\
&&+\rho(T(u_4), T(u_5))(\lambda(T(u_1), T(u_2))(u_3))+\rho(T(u_4), T(u_5))(\mu(T(u_1), T(u_3))(u_2))\nonumber\\
&&+\rho(T(u_4), T(u_5))(\rho(T(u_2), T(u_3))(u_1)\nonumber.
 \end{eqnarray}
In the same way, we get successively,
 \begin{eqnarray}
  [[u_1, u_4, u_5]_T, u_2, u_3]_T&=&\lambda(T(u_1), T(u_4), T(u_5), u_2)(u_3)+\mu(T(u_1), T(u_4), T(u_5), T(u_3))(u_2)\nonumber\\
&&+\rho(T(u_2), T(u_3))(\lambda(T(u_1), T(u_4))(u_5))+\rho(T(u_2), T(u_3))(\mu(T(u_1), T(u_5))(u_4))\nonumber\\
&&+\rho(T(u_2), T(u_3))(\rho(T(u_4), T(u_5))(u_1)\nonumber.
\end{eqnarray}
\begin{eqnarray}
[u_1, [u_2, u_4, u_5]_T, u_3]_T&=&\lambda(T(u_1), [T(u_2), T(u_4), T(u_5)])(u_3)+\mu(T(u_1), T(u_3))\lambda(T(u_2), T(u_4))(u_5)\nonumber\\
&&+\mu(T(u_1), T(u_3))\mu(T(u_2), T(u_5))(u_4)+\mu(T(u_1), T(u_3))\rho(T(u_4), T(u_5))(u_2)\nonumber\\
&&+\rho([T(u_2), T(u_4), T(u_5)], T(u_3))(u_1)\nonumber.
 \end{eqnarray}
\begin{eqnarray}
 [u_1, u_2, [u_3, u_4, u_5]_T]_T&=&\lambda(T(u_1), T(u_2))\lambda(T(u_3), T(u_4))(u_5)+\lambda(T(u_1), T(u_2))\mu(T(u_3), T(u_5))(u_4)\nonumber\\
&&+\lambda(T(u_1), T(u_2))\rho(T(u_4), T(u_5))(u_3)+\mu(T(u_1), [T(u_3), T(u_4), T(u_5)])(u_2)\nonumber\\
&&+\rho(T(u_2), [T(u_3), T(u_4), T(u_5)])(u_1).\nonumber
\end{eqnarray}
Using axioms in Definition \ref{dr}, we get the conclusion of the first part. The second part comes from the definitions of the 
$O$-operator $T$ and the bracket $[-, -, -]_T$.
\end{proof}

Taking $\rho(y, z)x=\mu(x, z)y=\lambda(x, y)z=[x, y, z]$, we obtain the below consquence.
\begin{corollary}(\cite{IS})
 Given a Rota-Baxter operator $R: L\rightarrow L$ of weight $\lambda=0$ on a ternary Leibniz algebra $L$, we can make $L$ into another ternary Leibniz 
algebra with the bracket
\begin{eqnarray}
 &&\{x, y, z\}=[R(x), R(y), z]+[R(x), y, R(z)]+[x, R(y), R(z)],
\end{eqnarray}
for any $x, y, z\in L$.
\end{corollary}

\section{Modules, representations and deformations}



\subsection{Compatible modules}

\begin{definition}
 Let $(L, [-, -, -])$ be a ternary Leibniz algebra and $M$ a bimodule on $L$. A trilinear map $\omega : L\times L\times L\rightarrow M$ is said
 to be a $2$-cocycle on $L$ with values in $M$ if the following equality holds
\begin{eqnarray}
 l_1(\omega(x, y, z), t, u)+\omega([x, y, z], t, u)&=&l_3(x, y, \omega(z, t, u))+l_2(x, \omega(y, t, u), z)
+l_1(\omega(x, t, u), y, z)\nonumber\\
&&+\omega(x, y, [z, t, u])+\omega(x, [y, t, u], z)\nonumber
+\omega([x, t, u], y, z),
\end{eqnarray}
for any $x, y, z, t, u\in L$.
\end{definition}

Observe that if $\omega : L\times L\times L\rightarrow L$ is a $2$-cocycle which is a ternary bracket, the cocycle condition is nothing but 
the compatibility condition (\ref{cc}).
\begin{theorem}
 Let $(L, [-, -, -])$ be a ternary Leibniz algebra, $(M, l_1, l_2, l_3)$ a bimodule on $L$ and $\omega : L\times L\times L\rightarrow M$ a $2$-cocycle on $L$ 
with values in $M$. Then $L\oplus M$ is a ternary algebra with the bracket
\begin{eqnarray}
 [x+m_1, y+m_2, z+m_3]_\omega=[x, y, z]+l_1(m_1, y, z)+l_2(x, m_2, z)+l_3(x, y, m_3)+\omega(x, y, z),\nonumber
\end{eqnarray}
for all $x, y, z\in L, m_i\in M, i=1,2,3$. Moreover, the brackets $[-, -, -]_0$ and $[-, -, -]_\omega$ are compatible, 
where $0: L\times L\times L\rightarrow M$
 is the null $2$-cocycle.
\end{theorem}
\begin{proof}
We have to prove the ternary Leibniz identity for the bracket $[-, -, -]_\omega$. Indeed, for any $x, y, z\in L, m_i\in M$, $i=1, 2, 3, 4, 5$,
\begin{eqnarray}
&&\hspace{2cm} [[x+m_1, y+m_2, z+m_3]_\omega,  t+m_4, u+m_5]_\omega\nonumber\\
&&= [[x, y, z]+l_1(m_1, y, z)+l_2(x, m_2, z)+l_3(x, y, m_3)+\omega(x, y, z),  t+m_4, u+m_5]_\omega\nonumber\\
&&=[[x, y, z], t, u] +l_1\Big(l_1(m_1, y, z)+l_2(x, m_2, z)+l_3(x, y, m_3)+\omega(x, y, z),   t, u\Big)\nonumber\\
&&\qquad+l_2([x, y, z], m_4, u)+l_3([x, y, z], t, m_5)+\omega([x, y, z], t, u)\nonumber\\
&&=[[x, y, z], t, u] +l_1(l_1(m_1, y, z),  t, u)+l_1(l_2(x, m_2, z), t, u)+l_1(l_3(x, y, m_3), t, u)\nonumber\\
&&\qquad+l_1(\omega(x, y, z), t, u)  +l_2([x, y, z], m_4, u)+l_3([x, y, z], t, m_5)+\omega([x, y, z], t, u)\nonumber.
\end{eqnarray}
In the same maner, 
\begin{eqnarray}
[[x+m_1, t+m_4, u+m_5]_\omega,  y+m_2, z+m_3]_\omega&=&
[[x, t, u], y, z] +l_1(l_1(m_1, t, u),  y, z)+l_1(l_2(x, m_2, u), y, z)\nonumber\\
&&+l_1(l_3(x, t, m_5), y, z)+l_1(\omega(x, t, u), y, z)  +l_2([x, t, u], m_2, z)\nonumber\\
&&+l_3([x, t, u], y, m_3)+\omega([x, t, u], y, z).\nonumber
\end{eqnarray}
\begin{eqnarray}
 [x+m_1, y+m_2, [z+m_3,  t+m_4, u+m_5]_\omega]_\omega&=&[x, y, [z, t, u]]+l_1(m_1, y, [z, t, u])+l_2(x, m_2, [z, t, u])\nonumber\\
&&+l_3(x, y, l_1(m_3, t, u))+l_3(x, y, l_2(z, m_4, u))+l_3(x, y, l_3(z, t, m_5))\nonumber\\
&&+\omega(x, y, [z, t, u])+l_3(x, y, \omega(z, t, u))\nonumber
\end{eqnarray}
\begin{eqnarray}
 [x+m_1, [y+m_2,  t+m_4, u+m_5]_\omega,  z+m_3]_\omega&=&[x, [y, t, u], z]+l_1(m_1, [y, t, u], z)+l_2(x, [m_2, t, u], z)\nonumber\\
&&+l_2(x, l_2(y, m_4, u), z)+l_2(x, l_3(y, t, m_5), z)+l_3(x, [y, t, u], m_3)\nonumber\\
&&+l_2(x, \omega(y, t, u), z)+\omega(x, [y, t, u], z)\nonumber
\end{eqnarray}
For the second part, we have, by using above result,
\begin{eqnarray}
 &&\qquad[[x+m_1, y+m_2, z+m_3],  t+m_4, u+m_5]_\omega+[[x+m_1, y+m_2, z+m_3]_\omega,  t+m_4, u+m_5]\nonumber\\
&&\qquad-[x+m_1, y+m_2, [z+m_3,  t+m_4, u+m_5]]_\omega-[x+m_1, y+m_2, [z+m_3,  t+m_4, u+m_5]_\omega]\nonumber\\
&&\qquad-[x+m_1, [y+m_2,  t+m_4, u+m_5],  z+m_3]_\omega-[x+m_1, [y+m_2,  t+m_4, u+m_5]_\omega,  z+m_3]\nonumber\\
&&\qquad-[[x+m_1, t+m_4, u+m_5],  y+m_2, z+m_3]_\omega-[[x+m_1, t+m_4, u+m_5]_\omega,  y+m_2, z+m_3]\nonumber
\end{eqnarray}
\begin{eqnarray}
&&=[[x, y, z], t, u] +l_1(l_1(m_1, y, z),  t, u)+l_1(l_2(x, m_2, z), t, u)+l_1(l_3(x, y, m_3) t, u)\nonumber\\
&&\quad \quad +l_2([x, y, z], m_4, u)+l_3([x, y, z], t, m_5)+\omega([x, y, z], t, u)\nonumber\\
&&\quad+[[x, y, z], t, u] +l_1(l_1(m_1, y, z),  t, u)+l_1(l_2(x, m_2, z), t, u)+l_1(l_3(x, y, m_3), t, u]\nonumber\\
&&\quad\quad+l_1(\omega(x, y, z), t, u)  +l_2([x, y, z], m_4, u)+l_3([x, y, z], t, m_5)\nonumber\\
&&\quad-[x, y, [z, t, u]]-l_1(m_1, y, [z, t, u])-l_2(x, m_2, [z, t, u])-l_3(x, y, l_1(m_3, t, u))\nonumber\\
&&\quad\quad-l_3(x, y, l_2(z, m_4, u))-l_3(x, y, l_3(z, t, m_5))-l_3(x, y, \omega(z, t, u))\nonumber\\
&&\quad-[x, y, [z, t, u]]-l_1(m_1, y, [z, t, u])-l_2(x, m_2, [z, t, u])-l_3(x, y, l_1(m_3, t, u))-l_3(x, y, l_2(z, m_4, u))\nonumber\\
&&\quad\quad-l_3(x, y, l_3(z, t, m_5))-\omega(x, y, [z, t, u])\nonumber\\
&&\quad-[x, [y, t, u], z]-l_1(m_1, [y, t, u], z)-l_2(x, l_1(m_2, t, u), z)-l_2(x, l_2(y, m_4, u), z)-l_2(x, l_3(y, t, m_5), z)\nonumber\\
&&\quad\quad-l_3(x, [y, t, u], m_3)-\omega(x, [y, t, u], z)\nonumber\\
&&\quad-[x, [y, t, u], z]-l_1(m_1, [y, t, u], z)-l_2(x, l_1(m_2, t, u), z)-l_2(x, l_2(y, m_4, u), z)-l_2(x, l_3(y, t, m_5), z)\nonumber\\
&&\quad\quad-l_3(x, [y, t, u], m_3)-l_2(x, \omega(y, t, u), z)\nonumber\\
&&\quad-[[x, t, u], y, z] -l_1(l_1(m_1, t, u),  y, z)-l_1(l_2(x, m_2, u), y, z)-l_1(l_3(x, t, m_5), y, z)\nonumber\\
&&\quad\quad-l_2([x, t, u], m_2, z)-l_3([x, t, u], y, m_3)-\omega([x, t, u], y, z)\nonumber\\
&&\quad-[[x, t, u], y, z] -l_1(l_1(m_1, t, u),  y, z)-l_1(l_2(x, m_2, u), y, z)-l_1(l_3(x, t, m_5), y, z)\nonumber\\
&&\quad\quad-l_1(\omega(x, t, u), y, z)-l_2([x, t, u], m_2, z)-l_3([x, t, u], y, m_3).\nonumber
\end{eqnarray}
The left hand side vanishes by ternary Leibniz identity, ternary Leibniz modules axioms and $2$-cocycle condition.
\end{proof} 

\begin{corollary}(\cite{IS})\label{ms}
 Let $(M, l_1, l_2, l_3)$ be a  bimodule over a ternary Leibniz algebra $L$. Then $L\oplus M$ is a ternary Leibniz algebra with respect to the 
multiplication 
\begin{eqnarray}
 [x_1+m_1, x_2+m_2, x_3+m_3]_0=[x_1, x_2, x_3]+l_1(m_1, x_2, x_3)+l_2(x_1, m_2, x_3)+l_3(x_1, x_2, m_3),\nonumber
\end{eqnarray}
for any $x_i+m_i\in L\oplus M, i=1, 2, 3$.
\end{corollary}
\begin{definition}
 Let $(L, [-, -, -]_1, [-, -, -]_2)$ be a compatible Leibniz algebra. A compatible bimodule over $L$ is a linear vector space $M$ together with six linear
maps $l_i^j, i=1, 2, 3, j=1, 2$ 
\begin{eqnarray}
\begin{array}{lll}
l_1^1  : M\times L\times L\rightarrow M, \: (m, x, y)\mapsto l_1(m, x, y) & ; &   l_1^2: M\times L\times L\rightarrow M, \: (m, x, y)\mapsto l_1^2(m, x, y)\\
 l_2^1 : L\times M\rightarrow L, \: (x, m, y)\mapsto l_2^1(x, m, y) & ; &  l_2^2 : L\times M\rightarrow L, \: (x, m, y)\mapsto l_2^2(x, m, y),\nonumber\\
l_3^1 : L\times L\rightarrow M, \: (x, y, m)\mapsto l_3^1(x, y, m) & ; &  l_3^2 : L\times L\rightarrow M, \: (x, y, m)\mapsto l_3^2(x, y, m),
\end{array}
\end{eqnarray}
such that $(M, l_1^1, l_2^1, l_3^1)$, $(M, l_1^2, l_2^2, l_3^2)$ are bimodules over $(L, [-, -, -]_1)$ and $(L, [-, -, -]_2)$ respectively, and
\begin{eqnarray}
 l_1^2(l_1^1(m, x, y), z, t)+l_1^1(l_1^2(m, x, y), z, t)&=&l_1^2(m, x, [y, z, t]_1)+l_1^1(m, x, [y, z, t]_2)+l_1^2(m, [x, z, t]_1, y)\nonumber\\
&&+l_1^1(m, [x, z, t]_2, y)+l_1^2(l_1^1(m, z, t), x, y)+l_1^1(l_1^2(m, z, t), x, y),\\
l_1^2(l_2^1(x, m, y), z, t)+l_1^1(l_2^2(x, m, y), z, t)&=&l_2^2(x, m, [y, z, t]_1)+l_2^1(x, m, [y, z, t]_2)+l_2^2(x, l_1^1(m, z, t), y)\nonumber\\
&&+l_2^1(x, l_1^2(m, z, t), y)+l_2^2([[x, z, t]_1, m, y)+l_2^1([x, z, t]_2, m, y),\\
l_1^2(l_3^1(x, y, m), z, t)+l_1^1(l_3^2(x, y, m), z, t)&=&l_3^2(x, y, l_1^1(m, z, t))+l_3^1(x, y, l_1^2(m, z, t))+l_3^2(x, [y, z, t]_1, m)\nonumber\\
&&+l_3^1(x, [y, z, t]_2, m)+l_3^2([x, z, t]_1, y, m)+l_3^1([[x, z, t]_2, y, m),\\
l_2^2([x, y, z]_1, m, t)+l_2^1([x, y, z]_2, m, t)&=&l_3^2(x, y, l_2^1(z, m, t))+l_3^1(x, y, l_2^2(z, m, t))+l_2^2(x, l_2^1(y, m, t), z)\nonumber\\
&&+l_2^1(x, l_2^2(y, m, t), z)+l_1^2(l_2^1(x, m, t), y, z)+l_1^1(l_2^2(x, m, t), y, z),\\
l_3^2([x, y, z]_1, t, m)+l_3^1([x, y, z]_2, t, m)&=&l_3^2(x, y, l_3^1(z, t, m))+l_3^1(x, y, l_3^2(z, t, m))+l_2^2(x, l_3^1(y, t, m), z)\nonumber\\
&&+l_2^1(x, l_3^2(y, t, m), z)+l_1^2(l_3^1(x, t, m), y, z)+l_1^1(l_3^2(x, t, m), y, z).
\end{eqnarray}
for any $x, y, z, t\in L$, $m\in M$,
\end{definition}

\begin{proposition}\label{pi}
 Let $(L, [-, -, -]_1, [-, -, -]_2)$ be a compatible ternary Leibniz algebra and $(M, l_1^i, l_2^i, l_3^i), i=1, 2$ be a compatible 
$L$-module. Then, the direct sum $L\oplus M$ has a compatible Leibniz algebra structure with the following multiplications :
\begin{eqnarray}
 [x+m_1, y+m_2, y+m_3]^i:=[x, y, z]_1+l_1^i(m_1, x, y)+l_2^i(x, m_2, y)+l_3^i(x, y, m_3), 
\end{eqnarray}
for any $x, y, z\in L, m_j\in M, j=1,2, 3, i=1, 2$.
\end{proposition}
\begin{proof}
 It follows from a direct computation.
\end{proof}

\subsection{Compatible representations}
Now we deal with the study of representations of ternary Leibniz  algebras.

\begin{proposition}\label{}
 Let $(L, [-, -, -])$ be a ternary Leibniz algebra and $(M, \lambda, \mu, \rho)$ be a representation of 
$L$ on $M$. Then, the direct sum $L\oplus M$ has a Leibniz algebra structure with the product :
\begin{eqnarray}
 [x+m_1, y+m_2, y+m_3]:=[x, y, z]+\rho(y, z)(m_1)+\mu(x, z)(m_2)+\lambda(x, y)(m_3), \nonumber
\end{eqnarray}
for any $x, y, z\in L, m_i, i=1,2, 3$.
\end{proposition}
\begin{proof}
The proof is the same as in that of Corollary \ref{ms}.
\end{proof}

Now we discuss on the compatibility of representations on compatible ternary Leibniz algebras. 

\begin{definition}
Let $(L, [-, -, -])$ be a ternary Leibniz algebra, $(V, \lambda_i, \mu_i, \rho_i), i=1, 2$ be two representations of $L$ on $V$. We say that
 the representations
 $(V, \lambda_1, \mu_1, \rho_1)$ and $(V, \lambda_2, \mu_2, \rho_2)$ are compatible if $(V, \lambda_1+\lambda_2, \mu_1+\mu_2, \rho_1+\rho_2)$
 is also a representation of $L$ on $V$. In other words, the representations
 $(V, \lambda_1, \mu_1, \rho_1)$ and $(V, \lambda_2, \mu_2, \rho_2)$ are compatible if, for any $x, y, z, t\in L$, we have 
\begin{eqnarray}
\rho_1(z, t)\lambda_2(x, y)+\rho_2(z, t)\lambda_1(x, y)-\lambda_1(x, y)\rho_2(z, t)-\lambda_2(x, y)\rho_1(z, t)&=&0\\
\rho_1(z, t)\mu_2(x, y)+\rho_2(z, t)\mu_1(x, y)-\mu_1(x, y)\rho_2(z, t)-\mu_2(x, y)\rho_1(z, t)&=&0\\
\rho_1(z, t)\rho_2(x, y)+\rho_2(z, t)\rho_1(x, y)-\rho_1(x, y)\rho_2(z, t)-\rho_2(x, y)\rho_1(z, t)&=&0,\\
 \lambda_1(x, y)\lambda_2(z, t)+\lambda_2(x, y)\lambda_1(z, t)+\mu_1(x, z)\lambda_2(y, t)+\mu_2(x, z)\lambda_1(y, t)
+\rho_1(y, z)\lambda_2(x, t)+\rho_2(y, z)\lambda_1(x, t)&=&0,\\
\lambda_1(x, y)\mu_2(z, t)+\lambda_2(x, y)\mu_1(z, t)+\mu_1(x, z)\mu_2(y, t)+\mu_2(x, z)\mu_1(y, t)+\rho_1(y, z)\mu_2(x, t).
+\rho^2(y, z)\mu_1(x, t)&=&0.
\end{eqnarray}
\end{definition}

\begin{theorem}
 Let $(L, [-, -, -]_1, [-, -, -]_2)$ be a compatible ternary Leibniz algebra and $(V, \lambda_i, \mu_i, \rho_i), i=1, 2$ be a compatible
 representation of $L$ on $V$. Then, the direct sum $L\oplus V$ carries a compatible ternary Leibniz structure with 
\begin{eqnarray}
 [x+m_1, y+m_2, z+m_3]_i=[x, y, z]_i+\lambda_i(x, y)(m_3)+\mu_i(x, z)(m_2)+\rho_i(y, z)(m_1), 
\end{eqnarray}
for any $x, y, z\in L, m_j\in V, i=1, 2, j=1, 2, 3$.
\end{theorem}
\begin{proof}
 It comes from direct calculation.
\end{proof}

\subsection{Infinitesimal deformation with a representation}

\begin{definition}
 Let $(L, [-, -, -])$ be a ternary algebra, $V$ a linear vector space, $\lambda, \mu, \rho : L\times L\rightarrow End(V)$ a representation of $L$ 
on $V$. Let $\omega^i : L\times L\rightarrow L$, $\omega_\lambda^i, \omega_\mu^i, \omega_\rho^i : L\times L\rightarrow End(V), i=1, 2$ be bilinear maps.
Consider a $s$-parametrized family bracket operations and linear maps
\begin{eqnarray}
 [x, y, z]_s&=&[x, y, z]+s\omega^1(x, y, z)+s^2\omega^2(x, y, z)\\
\lambda_s(x, y)&=&\lambda(x, y)+s\omega_\lambda^1(x, y, z)+s^2\omega_\lambda^2(x, y)\\
\mu_s(x, y)&=&\mu(x, y)+s\omega_\mu^1(x, y, z)+s^2\omega_\mu^2(x, y)\\
\rho_s(x, y)&=&\rho(x, y)+s\omega_\rho^1(x, y, z)+s^2\omega_\rho^2(x, y),
\end{eqnarray}
for any $x, y\in L$.
If $(L, [-, -, -]_s)$ are ternary Leibniz algebras and $(\lambda_s^i, \mu_s^i, \rho_s^i)$ are representations of $(L, [-, -, -]_s)$ on $V$, 
we say that $(\omega^i, \omega_\lambda^i, \omega_\mu^i, \omega_\rho^i)$ generates a one parameter infinitesimal deformation of the ternary 
Leibniz algebra $(L, [-, -, -])$ with the representation $(V, \lambda, \mu, \rho)$. 
\end{definition}
We denote a one-parameter infinitesimal deformation of a ternary Leibniz algebra $(L, [-, -, -])$ with the representation 
$(V, \lambda, \mu, \rho)$ by $(L, [-, -, -]_s, \lambda_s, \mu_s, \rho_s)$.

\begin{theorem}
Let $(L, [-, -, -])$ be a ternary Leibniz algebra. The quintuple $(L, [-, -, -]_s, \lambda_s, \mu_s, \rho_s)$ is a one-parameter infinitesimal deformation of a ternary
Leibniz algebra $L$ with the representation $(V, \lambda, \mu, \rho)$ if and only if the following set of equations holds :
\begin{eqnarray}
 &&[\omega^1(x, y, z), t, u]+\omega^1([x, y, z], t, u)
=[x, y, \omega^1(z, t, u)]+\omega^1(x, y [z, t, u])+[x, \omega^1(y, t, u), z]+\omega^1(x, [y, t, u], z)\nonumber\\
&&\hspace{6cm}+[\omega^1(x, t, u), y, z]+\omega^1([x, t, u], y, z),\label{d1}\\
&&\omega^1(\omega^1(x, y, z), t, u)-\omega^1(x, y, \omega^1(z, t, u))-\omega^1(x, \omega^1(y, t, u), z)-\omega^1(\omega^1(x, t, u),y, z)\nonumber\\
&&\hspace{3cm}=-[\omega^2(x, y, z), t, u]-\omega^2([x, y, z], t, u)+[x, y, \omega^2(z, t, u)]+\omega^2(x, y, [z, t, u])
+[x, \omega^2(y, t, u), z]\nonumber\\
&&\hspace{3,5cm}+\omega^2(x, [y, t, u], z)+[\omega^2(x, t, u), y, z]+\omega^2([x, t, u], y, z)\label{d2}\\
&&\omega^1(\omega^2(x, y, z), t, u)+\omega^2(\omega^1(x, y, z), t, u)=\omega^1(x, y, \omega^2(z, t, u))+\omega^2(x, y, \omega^1(z, t, u))
+\omega^1(x, \omega^2(y, t, u), z)\nonumber\\
&&\hspace{6cm}+\omega^2(x, \omega^1(y, t, u), z)+\omega^1(\omega^2(x, t, u),y, z)+\omega^2(\omega^1(x, t, u), y, z)\label{d3}\\
&&\omega^2(\omega^2(x, y, z), t, u)=\omega^2(x, y, \omega^2(z, t, u))+\omega^2(x, \omega^2(y, t, u), z)+\omega^2(\omega^2(x, t, u),y, z)\label{d3'}.
\end{eqnarray}
\begin{eqnarray}
&& \lambda(\omega^1(x, y, z), t)+\omega_\lambda^1([x, y, z], t)=\lambda(x, y)\omega_\lambda^1(z, t)+\omega_\lambda^1(x, y)\lambda(z, t)
+\mu(x, z)\omega_\lambda^1(y, t)+\omega_\mu^1(x, z)\lambda(y, t)\label{d4}\\
&&\hspace{6cm}+\omega_\rho^1(y, z)\lambda(x, t)+\rho(y, z)\omega_\lambda^1(x, t),\nonumber\\
&&\omega_\rho^1(z, t)\mu(x, y)+\rho(z, t)\omega_\mu^1(x, y)=\mu(x, \omega^1(y, z, t))+\omega_\mu^1(x, [y, z, t])
+\mu(x, y)\omega_\rho^1(z, t)+\omega_\mu^1(x, y)\rho(z, t)\\
&&\hspace{6cm}+\omega_\mu^1([x, z, t], y)+\mu(\omega^1(x, z, t), y),\nonumber\\
&&\rho(z, t)\omega_\lambda^1(x, y)+\omega_\rho^1(z, t)\lambda(x, y)=\lambda(x, y)\omega_\rho^1(z, t)+\omega_\lambda^1(x, y)\rho(z, t)
+\lambda(x, \omega^1(y, z, t))+\omega_\lambda^1(x, [y, z, t])\\
&&\hspace{6cm}+\lambda(\omega^1(x, z, t), y)+\omega_\lambda^1([x, z, t], y)\nonumber\\
&&\mu(\omega^1(x, y, z), t)+\omega_\mu^1([x, y, z], t)=\lambda(x, y)\omega_\mu^1(z, t)+\omega_\lambda^1(x, y)\mu(z, t)
+\mu(x, z)\omega_\mu^1(y, t)
+\omega_\mu^1(x, z)\mu(y, t)\\
&&\hspace{6cm}+\rho(x, z)\omega_\mu^1(x, t)+\omega_\rho^1(x, z)\mu(x, t)\nonumber\\
&&\rho(z, t)\omega_\rho^1(x, y)+\omega_\rho^1(z,t)\rho(x, y)=\rho(x, \omega^1(y, z, t))+\omega_\rho(x, [y, z, t])+\rho(\omega^1(x, z, t), y)
+\omega_\rho^1([x, z, t], t)\label{d5}\\
&&\hspace{6cm}+\rho(x, y)\omega_\rho^1(z, t)+\omega_\rho(x, y)\rho(z, t)\nonumber.
\end{eqnarray}
\begin{eqnarray}
&&\lambda(\omega(x, y, z), t)+\omega_\lambda^2([x, y, z], t)-\lambda(x, y)\omega_\lambda^2(z, t)
-\omega_\lambda^2(x, y)\lambda(z, t)-\mu(x, z)\omega_\lambda^2(y, t)-\omega_\mu^2(x, z)\lambda(y, t)\label{d6}\\
&&\quad-\rho(y, z)\omega_\lambda^2(x, t)-\omega_\rho^2(y, z)\lambda(x, t)
=\omega_\lambda^1(x, y)\omega_\lambda^1(z, t)+\omega_\mu^1(x, z)\omega_\lambda^1(y, t)
+\omega_\rho^1(y, z)\omega_\lambda^1(x, t)-\omega_\lambda^1(\omega(x, y, z), t),\nonumber\\
&&\rho(z, t)\omega_\mu^2(x, y)+\omega_\rho^2(z, t)\mu(x, y)-\mu(x, y)\omega_\rho^2(z, t)
-\mu(x, \omega^2(y, z, t))-\omega_\mu^2(x, [y, z, t])-\omega_\mu^2(x, y)\rho(z, t)\\
&&\qquad-\mu(\omega^2(x, z, t), y)-\omega_\mu^2([x, z, t], y)=\omega_\mu^1(x, \omega^1(y, z, t))+\omega_\mu^1(x, y)\omega_\rho^1(z, t)
+\omega_\mu^1(\omega^1(x, z, t), y)-\omega_\rho^1(z, t)\omega_\mu^1(x, y)\nonumber\\
&&\rho(z, t)\omega_\lambda^2(x, y)+\omega_\rho^2(z, t)\lambda(x, y)-\lambda(x, y)\omega_\rho^2(z, t)-\omega_\lambda^2(x, y)\rho(z, t)
-\lambda(x, \omega^2(y, z, t))-\omega_\lambda^2(x, [y, z, t])\\
&&\qquad-\lambda(\omega^2(x, z, t), y)-\omega_\lambda^2([x, z, t], y)=\omega_\lambda^1(x, \omega^1(y, z, t))+\omega_\lambda^1(\omega^1(x, z, t), y)
+\omega_\lambda^1(x, y)\omega_\rho^1(z, t)-\omega_\rho^1(z, t)\omega_\lambda^1(x, y)\nonumber\\
&&\mu(\omega^2(x, y, z), t)+\omega_\mu^2([x, y, z], t)-\lambda(x, y)\omega_\mu^2(z, t)-\omega_\lambda^2(x, y)\mu(z, t)-\mu(x, z)\omega_\mu^2(y, t)
-\omega_\mu^2(x, z)\mu(y, t)\\
&&\quad-\rho(x, z)\omega_\mu^2(x, t)-\omega_\rho^2(x, z)\mu(x, t)=\omega_\lambda^1(x, y)\omega_\mu^1(z, t)+\omega_\mu^1(x, z)\omega_\mu^1(y, t)+
\omega_\rho^1(x, z)\omega_\mu^1(x, t)-\omega_\mu^1(\omega^1(x, y, z), t)\nonumber\\
&&\rho(z, t)\omega_\rho^2(x, y)+\omega_\rho^2(z, t)\rho(x, y)-\rho(x, \omega^2(y, z, t))-\omega_\rho^2(x, [y, z, t])-\rho(\omega^2(x, z, t), y)
-\omega_\rho^2([x, z, t], y)\label{d7}\\
&&\quad-\rho(x, y)\omega_\rho^2(z, t)-\omega_\rho^2(x, y)\rho(z, t)=\omega_\rho^1(x, \omega^1(y, z, t))+\omega_\rho^1(\omega^1(x, z, t), y)
+\omega_\rho^1(x, y)\omega_\rho^1(z, t)-\omega_\rho^1(z, t)\omega_\rho^1(x, y)\label{d71}
\end{eqnarray}
\begin{eqnarray}
&&\omega_\lambda^1(\omega^2(x, y, z), t)+\omega_\lambda^2(\omega^1(x, y, z), t)=\omega_\lambda^1(x, y)\omega_\lambda^2(z, t)
+\omega_\lambda^2(x, y)\omega_\lambda^1(z, t)+\omega_\mu^1(x, z)\omega_\lambda^2(y, t)+\omega_\mu^2(x, z)\omega_\lambda^1(y, t)\nonumber\\
&&\hspace{6cm}+\omega_\rho^1(y, z)\omega_\lambda^2(x, t)+\omega_\rho^2(y, z)\omega_\lambda^1(x, t)\label{d8}\\
&&\omega_\rho^1(z, t)\omega_\mu^2(x, y)+\omega_\rho^2(z, t)\omega_\mu^1(x, y)=\omega_\mu^1(x, \omega(y, z, t))
+\omega_\mu^2(x, \omega^1(x, z, t))+\omega_\mu^1(x, y)\omega_\rho^2(z, t)+\omega_\mu^2(x, y)\omega_\rho^1(z, t)\nonumber\\
&&\hspace{6cm}+\omega_\mu^1(\omega^2(x, z, t), y)+\omega_\mu^2(\omega^1(x, z, t), y)\\
&&\omega_\rho^1(z, t)\omega_\lambda^2(x, y)+\omega_\rho^2(z, t)\omega_\lambda^1(x, y)
=\omega_\lambda^1(x, y)\omega_\rho^2(z, t)-\omega_\lambda^2(x, y)\omega_\rho^1(z, t)+\omega_\lambda^1(x, \omega^2(y, z, t))
+\omega_\lambda^2(x, \omega^1(y, z, t))\nonumber\\
&&\hspace{6cm}+\omega_\lambda^1(\omega^2(x, z, t), y)+\omega_\lambda^2(\omega^1(x, z, t), y)\\
&&\omega_\mu^1(\omega^2(x, y, z), t)+\omega_\mu^2(\omega^1(x, y, z), t)=\omega_\lambda^1(x, y)\omega_\mu^2(z, t)+\omega_\lambda^2(x, y)\omega_\mu^1(z, t)
+\omega_\mu^1(x, z)\omega_\mu^2(y, t)+\omega_\mu^2(x, z)\omega_\mu^1(y, t)\nonumber\\
&&\hspace{6cm}+\omega_\rho^1(x, z)\omega_\mu^2(x, t)+\omega_\rho^2(x, z)\omega_\mu^1(x, t)\\
&&\omega_\rho^1(z, t)\omega_\rho^2(x, y)+\omega_\rho^2(z, t)\omega_\rho^1(x, y)=\omega_\rho^1(x, \omega^2(y, z, t))
+\omega_\rho^2(x, \omega^1(y, z, t), y)+\omega_\rho^1(\omega^2(x, z, t), y)+\omega_\rho^2(\omega^1(x, z, t), y)\nonumber\\
&&\hspace{6cm}+\omega_\rho^1(x, y)\omega_\rho^2(z, t)+\omega_\rho^2(x, y)\omega_\rho^1(z, t)\label{d9}, 
\end{eqnarray}
for any $x, y, z, t, u\in L$.
\end{theorem}
\begin{proof}
 It follows from direct computation by writting the ternary Leibniz identity for $[-, -, -]_s$ and the representation axioms for
 $([-, -, -]_s, \lambda_s, \mu_s, \rho_s)$.
\end{proof}

We have the following observations.
\begin{remark}
\begin{enumerate}
 \item [1)]
Equality (\ref{d1}) means that $\omega^1$ is a $2$-cocycle on $(L, [-, -, -])$ with values in $L$. 
\item [2)]
The left hand side of (\ref{d2}) is the condition for $\omega^1$ to be a ternary Leibniz bracket while the right hand side is 
the compatibility condition of $[-, -, -]$ and $\omega^2$.
\item [3)] 
Equality (\ref{d3}),  is the compatibility condition of the brackets $\omega^1$ and $\omega^2$ or $\omega^1$ is a $2$-cocycle on $(L, \omega^2)$
with values in $L$. 
\item [4)]
Equality (\ref{d3'}) means that $\omega^2$ is a ternary Leibniz bracket.
\end{enumerate}
\end{remark}

\begin{corollary}
Suppose that $\omega^1$ is a ternary Leibniz bracket on $L$ and that $(L, [-, -, -]_s, \lambda_s, \mu_s, \rho_s)$ is
a one-parameter infinitesimal deformation of a ternary Leibniz algebra $(L, [-, -, -])$ with the representation 
$(V, \lambda, \mu, \rho)$. Then
\begin{enumerate}
 \item [1)]
the brackets $[-, -, -]$ and $\omega^1$ are compatible ternary brackets, 
\item [2)]
 the brackets $[-, -, -]$ and $\omega^2$ are compatible ternary brackets,
\item [3)] 
The brackets $\omega^1$ and $\omega^2$ are compatible ternary brackets,
\item [4)] 
the quadruple $(\lambda+\omega_\lambda^1, \mu+\omega_\mu^1, \rho+\omega_\rho^1)$ is a representation of 
$(L, [-, -, -]+\omega^1)$ on $V$.
\item [5)]
the quadruple $(\lambda+\omega_\lambda^2, \mu+\omega_\mu^2, \rho+\omega_\rho^2)$ is a representation of 
$(L, [-, -, -]+\omega^2)$ on $V$.
\item [6)]
the quadruple $(\omega_\lambda^1+\omega_\lambda^2, \omega_\mu^1+\omega_\mu^2, \omega_\rho^1+\omega_\rho^2)$ is a representation of 
$(L, \omega^1+\omega^2)$ on $V$.
\end{enumerate}
\end{corollary}
\begin{proof}
\begin{enumerate}
 \item [1)]
It comes from equality (\ref{d1}). 
\item [2)] 
It comes from equalities (\ref{d2}) and (\ref{d3'}).
\item [3)]
 It comes from equalities (\ref{d3}) and (\ref{d3'}).
\item [4)] 
 It comes from equalities(\ref{d4})-(\ref{d5}).
\item [5)]
 It comes from equalities (\ref{d6})-(\ref{d7}).
\item [6)]
 It comes from equalities (\ref{d8})-(\ref{d9}). 
\end{enumerate}
\end{proof}

\begin{corollary}
  The quintuple $(L, [-, -, -]_s, \lambda_s, \mu_s, \rho_s)$ is a 
linear deformation ($\omega^2=\omega_{\lambda}^2=\omega_{\mu}^2=\omega_{\rho}^2=0)$ of a ternary
Leibniz algebra $(L, [-, -, -])$ with the representation $(V, \lambda, \mu, \rho)$ if and only if :
\begin{enumerate}
 \item [1)]
$\omega^1$ is a  ternary Leibniz bracket.
 \item [2)]
 $[-, -, -]$ and $\omega^1$ are compatible ternary Leibniz brackets.
\item [3)]
 $(\omega_\lambda^1, \omega_\mu^1, \omega_\rho^1)$ is a representation of 
$(L, \omega^1)$ on $V$.
\item [4)] 
 $(\lambda+\omega_\lambda^1, \mu+\omega_\mu^1, \rho+\omega_\rho^1)$ is a representation of 
$(L, [-, -, -]+\omega^1)$ on $V$.
\end{enumerate}
\end{corollary}

Now we have the following definition.
\begin{definition}
i) Two one-parameter infinitesimal deformations $(L, [-, -, -]_s, \lambda_s, \mu_s, \rho_s)$ and $(L, [-, -, -]'_s, \lambda'_s, \mu'_s, \rho'_s)$
of a ternary Leibniz algebra $(L, [-, -, -])$ with the representation $(V, \lambda, \mu, \rho)$ are equivalent if there exists an isomorphism 
$(Id_L+sN, Id_V+sT)$ from $(L, [-, -, -]_s, \lambda_s, \mu_s, \rho_s)$ to $(L, [-, -, -]'_s, \lambda'_s, \mu'_s, \rho'_s)$ i.e. for any $x,, y\in L$,
\begin{eqnarray}
 (Id_L+sN)[x, y, z]'_s&=&[(Id_L+sN)(x), (Id_L+sN)(y), (Id_L+sN)(z)]_s,\\
 (Id_V+sT)\lambda'_s(x, y)&=&\lambda_s((Id_L+sN)(x), (Id_L+sN)(y))\circ (Id_V+sT),\\
(Id_V+sT)\mu'_s(x, y)&=&\mu_s((Id_L+sN)(x), (Id_L+sN)(y))\circ (Id_V+sT),\\
(Id_V+sT)\rho'_s(x, y)&=&\rho_s((Id_L+sN)(x), (Id_L+sN)(y))\circ (Id_V+sT),
\end{eqnarray}
ii) A one-parameter infinitesimal deformation of a Leibniz algebra $(L, [-, -, -])$ with representation $(V, \lambda, \mu, \rho)$ is said to be 
trivial if it is equivalent to $(L, [-, -, -], \lambda, \mu, \rho)$.
\end{definition}

\begin{definition}\label{dnr}
A pair $(N, T)$, where $N\in End(L)$ and $T\in End(T)$, is called a Nijenhuis pair on a ternary Leibniz algebra $(L, [-, -, -])$ with a 
representation $(V, \lambda, \mu, \rho)$ if $N$ is a Nijenhuis operator on the ternary Leibniz algebra $L$ and $T$ satisfies the following 
conditions hold : 
\begin{eqnarray}
\lambda(N(x), N(y))T &=&T\Big(\lambda(N(x), N(y))+\lambda(N(x), y)T+\lambda(x, N(y))T\Big)\nonumber\\
&&-T^2\Big(\lambda(x, y)T+\lambda(N(x), y)+\lambda(x, N(y))\Big)+T^3\lambda(x, y),\\
\mu(N(x), N(y))T &=&T\Big(\mu(N(x), N(y))+\mu(N(x), y)T+\mu(x, N(y))T\Big)\nonumber\\
&&-T^2\Big(\mu(x, y)T+\mu(N(x), y)+\mu(x, N(y))\Big)+T^3\mu(x, y),\\
\rho(N(x), N(y))T &=&T\rho(\rho(N(x), N(y))+\rho(N(x), y)T+\rho(x, N(y))T\Big)\nonumber\\
&&-T^2\Big(\rho(x, y)T+\rho(N(x), y)+\rho(x, N(y))\Big)+T^3\rho(x, y),
\end{eqnarray}
for any $x, y\in L$.
\end{definition}

\begin{theorem}
 Let $(N, S)$ be a Nijenhuis pair on a ternary Leibniz algebra $(L, [-, -, -])$ with a representation $(V, \lambda, \mu, \rho)$. 
Then, a deformation
of $(L, [-, -, -], \lambda, \mu, \rho)$ can be obtained by putting
\begin{eqnarray}
 \omega^1(x, y, z)&=&[N(x), y, z]+[x, N(y), z]+[x, y, N(z)],\\
\omega^2(x, y, z)&=&[N(x), N(y), z]+[N(x), y, N(z)]+[x, N(y), N(z)]-N([N(x), y, z]+[x, N(y), z]+[x, y, N(z)]),\\
\omega_\lambda^1(x, y)&=&\lambda(x, y)T+\lambda(N(x), y)+\lambda(x, N(y))-T\lambda(x, y),\\
\omega_\lambda^2(x, y)&=&\lambda(N(x), N(y))+\lambda(N(x), y)T+\lambda(x, N(y))T-T\omega_\lambda^1(x, y),\\
\omega_\mu^1(x, y)&=&\mu(x, y)T+\mu(N(x), y)+\mu(x, N(y))-T\mu(x, y),\\
\omega_\mu^2(x, y)&=&\mu(N(x), N(y))+\mu(N(x), y)T+\mu(x, N(y))T-T\omega_\mu^1(x, y),\\
\omega_\rho^1(x, y)&=&\rho(x, y)T+\rho(N(x), y)+\rho(x, N(y))-T\rho(x, y),\\
\omega_\rho^2(x, y)&=&\rho(N(x), N(y))+\rho(N(x), y)T+\rho(x, N(y))T-T\omega_\rho^1(x, y),
\end{eqnarray}
Furthermore, this deformation is trivial.
\end{theorem}
\begin{proof}
 It follows from a straightforward computation.
\end{proof}

\begin{definition}
A pair $(N, T)$, where $N\in End(L)$ and $T\in End(T)$, is dual-Nijenhuis pair on a ternary Leibniz algebra $(L, [-, -, -])$ with a 
representation $(V, \lambda, \mu, \rho)$ if $N$ is a Nijenhuis operator on the ternary Leibniz algebra $L$ and $T$ satisfies the following 
conditions : 
\begin{eqnarray}
 \lambda(x, y)T^3&=&T\lambda(N(x), N(y))-\Big(\lambda(N(x), N(y))+T\lambda(N(x), y)+T\lambda(x, N(y))\Big)T\nonumber\\
&&+\Big(T\lambda(x, y)+\lambda(N(x), y)+\lambda(x, N(y))\Big)T^2,\\
\mu(x, y)T^3&=&T\mu(N(x), N(y))-\Big(\mu(N(x), N(y))+T\mu(N(x), y)+T\mu(x, N(y))\Big)T\nonumber\\
&&+\Big(T\mu(x, y)+\mu(N(x), y)+\mu(x, N(y))\Big)T^2,\\
\rho(x, y)T^3&=&T\rho(N(x), N(y))-\Big(\rho(N(x), N(y))+T\rho(N(x), y)+T\rho(x, N(y))\Big)T\nonumber\\
&&+\Big(T\rho(x, y)+\rho(N(x), y)+\rho(x, N(y))\Big)T^2,
\end{eqnarray}
for any $x, y\in L$.
\end{definition}

Let $\lambda^*, \mu^*, \rho^* : L\times L\rightarrow End(V^*)$ be the dual representation of $\lambda, \mu, \rho$ defined by
\begin{eqnarray}
 \langle \lambda^*(x, y)\alpha, v\rangle=-\langle \alpha, \lambda(x, y)v\rangle,
\langle \mu^*(x, y)\alpha, v\rangle=-\langle \alpha, \mu(x, y)v\rangle,
\langle \rho^*(x, y)\alpha, v\rangle=-\langle \alpha, \rho(x, y)v\rangle,
 \forall x, y\in L, v\in V, \alpha\in V^*. \nonumber
\end{eqnarray}

\begin{proposition}
 $(N, T)$ is a Nijenhuis pair on a ternary Leibniz algebra $(L, [-, -, -])$ with a representation $(V, \lambda, \mu, \rho)$ if and only if
$(N, T^*)$ is a dual-Nijenhuis pair on the ternary Leibniz algebra $(L, [-, -, -])$ with the representation $(V^*, \lambda^*, \mu^*, \rho^*)$.
\end{proposition}
\begin{proof}
 It comes from a direct computation by transposition.
\end{proof}

\begin{proposition}
 Let $(N, T)$ be a Nijenhuis on a ternary Leibniz algebra $(L, [-, -, -])$ with a representation $(V, \lambda, \mu, \rho)$. Then,
$N+T$ is a Nijenhuis operator on the ternary Leibniz algebra $L\oplus V$ with the multiplication is defined by
\begin{eqnarray}
 [x+m_1, y+m_2, y+m_3]:=[x, y, z]+\rho(y, z)(m_1)+\mu(x, z)(m_2)+\lambda(x, y)(m_3), \nonumber
\end{eqnarray}
for all $x, y, z\in L, m_i, i=1,2, 3$.
\end{proposition}
\begin{proof}
For any $x, y, z\in L, m, n, p\in V$ and by Definition \ref{nod} and Definition \ref{dnr}, we have :
\begin{eqnarray}
 &&\qquad\{(N+T)(x+m), (N+T)(y+n), (N+T)(z+p)\}=\{(N(x)+T(m), N(y)+T(n), N(z)+T(p)\}\nonumber\\
&&=[N(x), N(y), N(z)]+\rho(N(y), N(z))(T(m))+\mu(N(x), N(z))(T(n))+\lambda(N(x), N(y))(T(p))\nonumber\\
&&=N\Big([N(x), N(y), z]+[x, N(y), N(z)]+[N(x), y, N(z)]\Big)-N^2\Big([N(x), y, z]+[x, N(y), z]+[x, y, N(z)]\Big)+N^3([x, y, z])\nonumber\\
&&\quad+T\Big(\rho(N(y), N(z))(m)+\rho(N(y), z)(T(m))+\rho(y, N(z))(T(m))\Big)-T^2\Big(\rho(y, z)(T(m))+\rho(N(y), z)
+\rho(y, N(z))\Big)\nonumber\\
&&\quad+T^3\rho(y, z)(m)+T\Big(\mu(N(x), N(z))(n)+\mu(N(x), z)(T(n))+\mu(x, N(z))(T(n))\Big)\nonumber\\
&&\quad-T^2\Big(\mu(x, z)(T(n))+\mu(N(x), z)(n)+\mu(x, N(z))(n)\Big)+T^3\mu(x, z)(n)\nonumber\\
&&\quad+T\Big(\lambda(N(x), N(y))(p)+\lambda(N(x), y)T(p)+\lambda(x, N(y))(T(p))\Big)\nonumber\\
&&\quad-T^2\Big(\lambda(x, y)(T(p))+\lambda(N(x), y)(p)+\lambda(x, N(y))(p)\Big)+T^3\lambda(x, y)(p)\nonumber.
\end{eqnarray}
Next,
\begin{eqnarray}
 &&\qquad(N+T)\Big(\{(N+T)(x+m), (N+T)(y+n), z+p\}+\{(N+T)(x+m), y+n, (N+T)(z+p)\}\nonumber\\
&&\hspace{4cm}+\{x+m, (N+T)(y+n), (N+T)(z+p)\}\Big)\nonumber\\
&&=(N+T)\Big(\{(N(x)+T(m), N(y)+T(n), z+p\}+\{(N(x)+T(m), y+n, N(z)+T(p)\}\nonumber\\
&&\qquad+\{x+m, N(y)+T(n), N(z)+T(p)\}\Big)\nonumber\\
&&=(N+T)\Big([N(x), N(y), z]+\rho(N(y), z)(T(m))+\mu(N(x), z)(T(n))+\lambda(N(x), N(y))(p)\nonumber\\
&&\hspace{2cm}+[N(x), y, N(z)]+\rho(y, N(z))(T(m))+\mu(N(x), N(z))(n)+\lambda(N(x), y)(T(p))\nonumber\\
&&\hspace{2cm}+[x, N(y), N(z)]+\rho(N(y), N(z))(m)+\mu(x, N(z))(T(n))+\lambda(x, N(y))(T(p))\Big)\nonumber
\end{eqnarray}
\begin{eqnarray}
&&=N\Big([N(x), N(y), z]+[N(x), y, N(z)]+[x, N(y), N(z)]\Big)\nonumber\\
&&+T\Big(\lambda(N(x), N(y))(p)+\lambda(N(x), y)(T(p))+\lambda(x, N(y))(T(p))\Big)\nonumber\\
&&+T\Big(\mu(N(x), z)(T(n))+  \mu(N(x), N(z))(n) +\mu(x, N(z))(T(n))\Big)\nonumber\\
&&+T\Big(\rho(N(y), z)(T(m))+\rho(y, N(z))(T(m))++\rho(N(y), N(z))(m)\Big)\nonumber.
\end{eqnarray}
Then,
\begin{eqnarray}
&& (N+T)^2\Big(\{(N+T)(x+m), y+n, z+p\}+\{x+m, (N+T)(y+n), z+p\}+\{x+m, y+n, (N+T)(z+p)\}\Big)=\nonumber\\
&&=N^2\Big([N(x), y, z]+[x, N(y), z]+[x, y, N(z)]\Big)+T^2\Big(\lambda(N(x), y)(p)+\lambda(x, N(y))(T(p))+\lambda(x, y)(T(p))\Big)\nonumber\\
&&\qquad+T^2\Big(\mu(N(x), z)(n)+  \mu(x, z)(T(n)) +\mu(x, N(z))(n)\Big)
+T^2\Big(\rho(y, z)(T(m))+\rho(N(y), z)(m)+\rho(y, N(z))(m)\Big)\nonumber.
\end{eqnarray}
and
\begin{eqnarray}
 (N+T)^3(\{x+m, y+n, z+p\})=N^3([x, y, z])+T^3\Big(\lambda(x, y)(p)+\mu(x, z)(n)+\rho(y, z)(m)\Big)\nonumber.
\end{eqnarray}
It follows that 
\begin{eqnarray}
&&\hspace{3cm}\{(N+T)(x+m), (N+T)(y+n), (N+T)(z+p)\}=\nonumber\\
&&=(N+T)\Big(\{(N+T)(x+m), (N+T)(y+n), z+p\}+\{(N+T)(x+m), y+n, (N+T)(z+p)\}\nonumber\\
&&\quad+\{x+m, (N+T)(y+n), (N+T)(z+p)\}\Big)-(N+T)^2\Big(\{(N+T)(x+m), y+n, z+p\}\nonumber\\
&&\quad+\{x+m, (N+T)(y+n), z+p\}+\{x+m, y+n, (N+T)(z+p)\}\Big)+(N+T)^3(\{x+m, y+n, z+p\})\nonumber.
\end{eqnarray}
This ends the proof.
\end{proof}

%
%
\end{document}